\documentclass[11pt]{article}
\usepackage[utf8]{inputenc}
\usepackage{amsfonts}
\usepackage{mathrsfs}
\usepackage{bbm}
\usepackage{amsfonts}
\usepackage{amssymb,amsmath,graphicx}
\usepackage{float} \usepackage[colorlinks=true]{hyperref}
\hypersetup{urlcolor=blue, citecolor=red}

\allowdisplaybreaks

\begin{document}

\title{Initial boundary value problem of a class of pseudo-parabolic Kirchhoff equations with logarithmic nonlinearity
}
\author{Qiuting Zhao{\thanks{Corresponding author. E-mail: zhaoqt72@163.com}}
\\
{\small School of Mathematical Sciences, Dalian University of Technology, Dalian 116024, China} }
\maketitle
\newtheorem{theorem}{Theorem}
\newtheorem{definition}{Definition}[section]
\newtheorem{lemma}{Lemma}[section]
\newtheorem{proposition}{Proposition}[section]
\newtheorem{corollary}{Corollary}[section]
\newtheorem{remark}{Remark}
\renewcommand{\theequation}{\thesection.\arabic{equation}}
\catcode`@=11 \@addtoreset{equation}{section} \catcode`@=12
\maketitle{}
{\bf Abstract}\quad
In this paper, we consider the initial boundary value problem for a pseudo-parabolic Kirchhoff equation
with logarithmic nonlinearity. We use the potential well method to give a threshold result
of global existence and finite-time blow-up for the weak solutions with initial energy $J(u_0)\leq d$.
When the initial energy $J(u_0)>d$, we find another criterion for the vanishing solution and blow-up solution.
We also get the exponential decay rate of the global solution and life span of the blow-up solution.
Meanwhile, we study the corresponding stationary problem and establish a convergence relationship
between its ground state solution and the global solution.



\section{Introduction}

In this paper, we are concerned with the following initial boundary value problem
\begin{eqnarray}\label{1.1}
\left\{
    \begin{array}{llll}
        \displaystyle u_{t}-k\Delta u_{t}-M(\|\nabla u\|_p^p)\Delta_p u=|u|^{q-1}u\log|u|,
        &&(x,t)\in\Omega\times(0,T),
        \\
        \displaystyle u(x,t)=0, &&(x,t)\in\partial\Omega\times(0,T),
        \\
        \displaystyle u(x,0)=u_0(x), && x\in\Omega,
    \end{array}
\right.
\end{eqnarray}
where $\Omega\subset\Bbb{R}^n (n\geq1)$ is a bounded domain with smooth boundary,
$M(s)=a+bs$, $a>0$, $b>0$, $k=0$ or $1$, $u(x, t) : \Omega\times(0,T)\rightarrow\Bbb{R}$,
$u_0(x)\in W_{0}^{1,p}(\Omega)$ with $u_0(x)\neq0$.
The parameters $p$ and $q$ satisfy $2\leq p<2^*$ and $2p-1<q<p^*-2$, where
$p^*=\frac{np}{n-p}$ is the Sobolev conjugate of $p$.

Problem \eqref{1.1} belongs to pseudo-parabolic Kirchhoff equation \cite{Showalter1}
which has had a high profile in the study of population dynamics
in recent years \cite{Tuan}, where the diffusion coefficient $M(\cdot)$
expresses the dependence on the global population density in the environment
instead of the density at a local location.
Equations like \eqref{1.1} can also be used to model many other phenomena such as
nonlinear elasticity, non-stationary fluid, image recovery,...
(see \cite{APS,CHM,CSH} and the references therein).
The logarithmic nonlinear source term is also used to describe
the mathematical and physical phenomena such as Gravity-mediated
super symmetric fracture model, wave equations,...(see
\cite{TMF,ZLH,CaoEJDE2,CTJDE,JCJMAA,Zhoujun2} and the references therein).

In 1883, Kirchhoff first proposed the Kirchhoff equation with
Dirichlet boundary \cite{GKiff}. Then, the convex method \cite{HAL1,HAL2}
and the potential well method \cite{DHS,DHS0} are successively proposed.
In 1978, Lions introduced a functional analysis framework \cite{Lions}
for equations with Kirchhoff terms. Since then, more and more scholars
have paid attention to this field and done excellent works.
Han et al. \cite{Han2} studied the following Kirchhoff equation
\begin{align}\label{qy1}
u_{t}-M(\|\nabla u\|_2^2)\Delta u=|u|^{q-1}u.
\end{align}
They considered the global existence, uniqueness, finite time blow-up
and asymptotic behavior of solutions with subcritical, critical
and supercritical initial energy. The upper and lower bounds of
the blow-up time of the solutions were supplemented in \cite{Han4}.
As an extension of \eqref{qy1}, Jian Li et al. \cite{HanMMA} replaced $\Delta u$
in \eqref{qy1} with $\Delta_p u$, describing the impact of the p-Laplace operator.
By the potential well method, R.Z. Xu et al. \cite{XS} studied
the following pseudo-parabolic equation
\begin{align*}
u_{t}-\Delta u_t-\Delta u=u^p.
\end{align*}
They proved the invariance of some sets, global existence, blow-up and asymptotic behavior
of solutions with different initial energies. Di et al. \cite{DiAML} studied
the pseudo-parabolic equation
\begin{align}\label{qy3}
u_{t}-\upsilon\Delta u_t-div(|\nabla u|^{p(x)-2}\nabla u)=|u|^{q(x)-2}u
\end{align}
with variable exponents, provided sufficient conditions for blow-up solutions
and estimated the upper bound of blow-up time.
Liao \cite{LiaoM} proved that the solution is not global in time
when the initial energy is positive, which extended and improved
the results in \cite{DiAML}.
Zhu X et al. \cite{ZhuAML} further supplemented the study of
the solutions of \eqref{qy3} with high initial energy.
In \cite{Zhoujun2}, the authors considered the case where $\upsilon=1$
and both p and q in \eqref{qy3} are constants. Under different values of p, q
and initial energies, there are global existence solutions, blow-up solutions
and extinct solutions. At the same time, the lower bound estimation of
the growth rate of the infinite time blow-up solutions is given.
Chen and Tian \cite{CT} studied the initial boundary value problem of
the following semi-linear pseudo-parabolic equation
\begin{align}\label{ct1}
u_{t}-\Delta u_t-\Delta u=u\log u,
\end{align}
using the Logarithmic Sobolev inequality and potential well family to obtain
the existence, blow-up and isolate vacuum of the global solution of the equation.
Meanwhile, the authors and discuss the asympototic behavior of solution.
On the basis of \cite{CT}, Nhan and Truong \cite{NT} studied the following equation
\begin{align}\label{nt1}
u_{t}-\Delta u_t-\Delta_p u=|u|^{p-2}u\log u,
\end{align}
and gave the sufficient condition for the existence of global decay solutions and finite-time blow-up solutions.

From many previous works about the IBVP of Kirchhoff equations, pseudo-parabolic equations
and other parabolic equations, the global solutions converge to 0 as t tends to $\infty$
when the initial data satisfies some special conditions.
Recently, some mathematicians \cite{Zhoujun1,Han3} make further discussion on the asymptotic behavior of general global solutions
and find that it is relevant to the ground state solutions of the stationary problem.

Motivated by these words, the rest of this paper is arranged as follows. Section 2 introduces preliminary knowledge. Sections 3-4 and 6 introduce
the global existence and the finite blow-up of the solutions of \eqref{1.1} in the case of $J(u_0)<d$, $J(u_0)=d$ and $J(u_0)>d$.
Section 5 estimates the asymptotic behavior of the global solutions and the life span
of finite time blow-up solutions. The last section, Section 7, introduces the boundary value problem
and establishes the convergence relation between the global solutions of \eqref{1.1}
and the ground state solutions of \eqref{1.2}.

\section{Preliminary}

\qquad In this paper, we use the expressions $\|\cdot\|_p=\|\cdot\|_{L^p(\Omega)}$,
$\|\cdot\|_{W_0^{1,p}}=\|\cdot\|_{W_0^{1,p}(\Omega)}$
and $(u,v)=\int_{\Omega}u(x)v(x)dx$.
\begin{definition}\label{def2.1}
A function $u(x,t)$ is called a weak solution to \eqref{1.1} on $\Omega\times[0,T)$,
if $u(x,0)=u_0(x)\in W_0^{1,p}(\Omega)$, $u\in L^\infty(0,T;W_{0}^{1,p}(\Omega))$ with
$u_t\in L^2(0,T;L^2(\Omega))$ and satisfies
\begin{align*}
   (u_t, \varphi)+k(\nabla u_t, \nabla\varphi)
   +M(\|\nabla u\|_p^p)(|\nabla u|^{p-2}\nabla u, \nabla\varphi)
   =(|u|^{q-1}u\log|u|, \varphi),
\end{align*}
for any $\varphi\in W_0^{1,p}(\Omega)$.
\end{definition}
\qquad Using the potential well theory \cite{DHS,Liu,XS}, we introduce the potential energy functional
\begin{align}\label{1}
J(u)={\frac{a}{p}}\|\nabla u\|_p^p+{\frac{b}{2p}}\|\nabla u\|_p^{2p}
-{\frac{1}{q+1}}\int_\Omega{|u|^{q+1}\log|u|}dx+\frac{1}{(q+1)^2}\|u\|_{q+1}^{q+1},
\end{align}
and the Nehari functional
\begin{align}\label{2}
I(u)=a\|\nabla u\|_p^p+b\|\nabla u\|_p^{2p}-\int_\Omega{|u|^{q+1}\log|u|}dx.
\end{align}
\eqref{1} and \eqref{2} imply that
\begin{align}
\label{3}
&J(u)={\frac{1}{q+1}}I(u)+\left({\frac{a}{p}}-{\frac{a}{q+1}}\right)\|\nabla u\|_p^p
+\left({\frac{b}{2p}}-{\frac{b}{q+1}}\right)\|\nabla u\|_p^{2p}+\frac{1}{(q+1)^2}\|u\|_{q+1}^{q+1},
\\
\label{cao2}
&\frac{d}{dt}J(u)=-\|u_t\|_2^2-k\|\nabla u_t\|_2^2.
\end{align}
For any $\delta>0$, the modified Nehari functional can be defined as
\begin{align*}
I_{\delta}(u)=\delta(a\|\nabla u\|_p^p+b\|\nabla u\|_p^{2p})-\int_\Omega{|u|^{q+1}\log|u|}dx.
\end{align*}
Then we can define the Nehari manifold and the potential wells
\begin{align*}
&\mathscr{N}=\{u\in W_{0}^{1,p}(\Omega):I(u)=0, \|\nabla u\|_p\neq0\},
\\
&W=\{u\in W_{0}^{1,p}(\Omega):J(u)<d, I(u)>0\}\bigcup\{0\},
\\
&V=\{u\in W_{0}^{1,p}(\Omega):J(u)<d, I(u)<0\},
\\
&\mathscr{N}_{\delta}=\{u\in W_{0}^{1,p}(\Omega):I_{\delta}(u)=0, \|\nabla u\|_p\neq0\},
\\
&W_{\delta}=\{u\in W_{0}^{1,p}(\Omega):J(u)<d(\delta), I_{\delta}(u)>0\}\bigcup\{0\},
\\
&V_{\delta}=\{u\in W_{0}^{1,p}(\Omega):J(u)<d(\delta), I_{\delta}(u)<0\},
\end{align*}
where $d(\delta)$ is the depth of the potential well $W_{\delta}$ and
\begin{equation}\label{cao12}
d=d(1)=\inf\{J(u):u\in\mathscr{N}\},\qquad d(\delta)=\inf\{J(u):u\in\mathscr{N}_{\delta}\}.
\end{equation}

\begin{lemma}\label{lem1}
For any $u\in W_{0}^{1,p}(\Omega)$ with $\|\nabla u\|_p\neq0$, there hold
\\
{\rm(i)} $\lim\limits_{{\lambda\rightarrow 0}}J(\lambda u)=0$,
$\lim\limits_{{\lambda\rightarrow+\infty}}J(\lambda u)=-\infty$.
\\
{\rm(ii)} There exists a unique $\lambda^*>0$ such that
$\frac{d}{d\lambda}J(\lambda u)|_{\lambda=\lambda^*}=0$,
namely $\lambda^* u\in \mathscr{N}$.
Furthermore, $\frac{d}{d\lambda}J(\lambda u)|_{\lambda=\lambda^*}>0$ on $(0,\lambda^*)$,
$\frac{d}{d\lambda}J(\lambda u)|_{\lambda=\lambda^*}<0$ on $(\lambda^*,\infty)$,
namely $J(\lambda u)$ takes the maximum at $\lambda=\lambda^*$.
\end{lemma}
{\bf Proof}
(i) For any $u\in W_0^{1,p}(\Omega)$ and $\lambda>0$,
\begin{align}\label{4}
J(\lambda u)&=\lambda^p{\frac{a}{p}}\|\nabla u\|_p^p
+\lambda^{2p}{\frac{b}{2p}}\|\nabla u\|_p^{2p}
+{\frac{\lambda^{q+1}}{(q+1)^2}}\|u\|_{q+1}^{q+1}
-{\frac{1}{q+1}}\int_\Omega{|\lambda u|^{q+1}\log|\lambda u|}dx.
\end{align}
Since $q+1>2p$, thus $\lim\limits_{{\lambda\rightarrow 0}}J(\lambda u)=0$,
$\lim\limits_{{\lambda\rightarrow+\infty}}J(\lambda u)=-\infty$.

(ii) Derivative $J(\lambda u)$ with respect to $\lambda$, we have
\begin{align*}
\frac{d}{d\lambda}J(\lambda u)&=\lambda^{p-1}a\|\nabla u\|_p^p
+\lambda^{2p-1}b\|\nabla u\|_p^{2p}
-\int_\Omega{\lambda^q|u|^{q+1}\log|\lambda u|}dx
\\
&=\lambda^{q}\left(\frac{a}{\lambda^{q+1-p}}\|\nabla u\|_p^p
+\frac{b}{\lambda^{q+1-2p}}\|\nabla u\|_p^{2p}
-\int_\Omega{|u|^{q+1}\log|\lambda u|}dx\right).
\end{align*}
Set $g(\lambda)=\frac{a}{\lambda^{q+1-p}}\|\nabla u\|_p^p
+\frac{b}{\lambda^{q+1-2p}}\|\nabla u\|_p^{2p}
-\int_\Omega{|u|^{q+1}\log|\lambda u|}dx$, then
\begin{align*}
&\lim\limits_{{\lambda\rightarrow 0}}g(\lambda)=+\infty,
\lim\limits_{{\lambda\rightarrow+\infty}}g(\lambda)<0,
\\
&g'(\lambda)=-\frac{a(q+1-p)}{\lambda^{q+2-p}}\|\nabla u\|_p^p
-{\frac{b(q+1-2p)}{\lambda^{q+2-2p}}}\|\nabla u\|_p^{2p}
-\lambda^{-1}\|u\|_{q+1}^{q+1}<0.
\end{align*}
Therefore there exists a unique $\lambda^*>0$ such that $g(\lambda^*)=0$,
namely $\frac{d}{d\lambda}J(\lambda u)|_{\lambda=\lambda^*}=0$.
It is easily to find that $J(\lambda u)$ is strictly increasing on $(0,\lambda^*]$,
strictly decreasing on $(\lambda^*,\infty)$ and takes the maximum at $\lambda=\lambda^*$.
$\hfill\Box$

\begin{lemma}\label{lem2}
For any $u\in W_{0}^{1,p}(\Omega)$ with $\|\nabla u\|_p\neq0$,
$r(\delta)=(\frac{\delta a}{S^{q+2}})^{\frac{1}{q+2-p}}$,
where $S$ is the embedding coefficient of the Sobolev inequality $\|u\|_{q+2}\leq S\|\nabla u\|_p$, there hold
\\
{\rm(i)} If $0<\|\nabla u\|_p<r(\delta)$, then $I_\delta(u)>0$.
\\
{\rm(ii)} If $I_\delta(u)<0$, then $\|\nabla u\|_p>r(\delta)$.
\\
{\rm(iii)} If $I_\delta(u)=0$, then $\|\nabla u\|_p=0$ or $\|\nabla u\|_p>r(\delta)$.
\end{lemma}
{\bf Proof} (i)The Sobolev embedding inequality and $0<\|\nabla u\|_p<r(\delta)$ indicate that
\begin{align*}
\int_\Omega{|u|^{q+1}\log|u|}dx\leq\int_\Omega{|u|^{q+2}}dx
<S^{q+2}r^{q+2-p}(\delta)\|\nabla u\|^p_p=a\delta\|\nabla u\|^p_p
<a\delta\|\nabla u\|^p_p+b\delta\|\nabla u\|^{2p}_p,
\end{align*}
which means $I_\delta(u)>0$.

(ii) can be directly derived from (i).

(iii) If $\|\nabla u\|_p=0$, then $I_\delta(u)=0$. If $I_\delta(u)=0$ and $\|\nabla u\|_p\neq 0$,
then $a\delta\|\nabla u\|_p^p<\int_\Omega{|u|^{q+1}\log|u|}dx\leq S^{q+2}\|\nabla u\|_p^{q+2}$,
namely $\|\nabla u\|_p>r(\delta)$.
$\hfill\Box$

\begin{lemma}\label{lem3}
$d(\delta)$ satisfies
\\
{\rm(i)} $\lim\limits_{{\delta\rightarrow0^+}}d(\delta)=0$,
$\lim\limits_{{\delta\rightarrow+\infty}}d(\delta)=-\infty$.
\\
{\rm(ii)} $d(\delta)$ is monotonically increased on $0<\delta\leq 1$,
monotonically decreased on $\delta>1$ and the maximum is obtained at $\delta=1$.
\end{lemma}
{\bf Proof}
(i) For any $\lambda u\in \mathscr{N}_\delta$, we have
$$
\delta a\|\nabla u\|^p_p+\lambda^{p}\delta b\|\nabla u\|^{2p}_p
=\lambda^{q+1-p}\int_\Omega{\left(|u|^{q+1}\log|\lambda u|\right)}dx,
$$
which indicates
\begin{align}\label{5}
\delta={\frac{\lambda^{q+1-p}\int_\Omega{\left(|u|^{q+1}\log|\lambda u|\right)}dx}
{a\|\nabla u\|^p_p+b\lambda^p\|\nabla u\|^{2p}_p}}.
\end{align}
A directly computation on \eqref{5} show that $\lambda$ increases as $\delta$ increases,
$\delta$ increases as $\lambda$ increases and $\lim\limits_{{\delta\rightarrow {0^+}}}
\lambda(\delta)=0$, $\lim\limits_{{\delta\rightarrow +\infty}}\lambda(\delta)=+\infty$.
Thus from the definition of $d(\delta)$ and Lemma \ref{lem1}, we can get
\begin{align*}
   &0\leq\lim\limits_{{\delta\rightarrow {0^+}}}d(\delta)
   \leq\lim\limits_{{\delta\rightarrow {0^+}}}J(\lambda u)
   =\lim\limits_{{\lambda\rightarrow {0^+}}}J(\lambda u)=0,
   \\
   &\lim\limits_{{\delta\rightarrow +\infty}}d(\delta)
   \leq\lim\limits_{{\delta\rightarrow +\infty}}J(\lambda u)
   =\lim\limits_{{\lambda\rightarrow +\infty}}J(\lambda u)=-\infty.
\end{align*}
Therefore $\lim\limits_{{\delta\rightarrow {0^+}}}d(\delta)=0$ and
$\lim\limits_{{\delta\rightarrow +\infty}}d(\delta)=-\infty$.

(ii) Assume $0<\delta'<\delta''\leq1$ or $1<\delta''<\delta'$.
Let $h(\lambda)=J(\lambda(\delta)u)$ with $\lambda(\delta)u\in\mathscr{N}_{\delta}$, then
\begin{equation*}
   h'(\lambda)=\lambda^{p-1} a\|\nabla u\|_p^p+\lambda^{2p-1} b\|\nabla u\|_p^{2p}
   -\int_\Omega{\lambda^q|u|^{q+1}\log|\lambda u|}dx
   =\frac{1}{\lambda}I(\lambda u).
\end{equation*}
For any $u\in\mathscr{N}_{\delta''}$ with $\lambda(\delta'')=1$, set
$v=\lambda(\delta')u\in\mathscr{N}_{\delta'}$.
If $0<\delta'<\delta''\leq 1$, since $\lambda(\delta)$ increases as
$\delta$ increases, then
\begin{align*}
   J(u)-J(v)&=h(1)-h(\lambda(\delta'))=\frac{1-\lambda(\delta')}{\lambda^*}I(\lambda^*u)
   \\
   &=\frac{1-\lambda(\delta')}{\lambda^*}\left[a(1-\delta^*)\|\lambda^*\nabla u\|_p^p
   +b(1-\delta^*)\|\lambda^*\nabla u\|_p^{2p}\right]
   \\
   &>0,
\end{align*}
where $\lambda^*=\theta\lambda(\delta')+(1-\theta)\lambda(\delta'')$,
$\theta\in(0,1)$ and $\lambda^* u\in\mathscr{N}_{\delta^*}$.
\\
Therefore, for any $u\in\mathscr{N}_{\delta''}$, there exists
$v\in\mathscr{N}_{\delta'}$ such that $J(u)>J(v)$,
which leads to $d(\delta'')>d(\delta')$.
The case for $1<\delta''<\delta'$ is similarly.
$\hfill\Box$

\begin{lemma}\label{lem4}
For any $u\in W_{0}^{1,p}(\Omega)$ with $0<J(u)<d$,
the sign of $I_\delta(u)$ doesn't change for $\delta_1<\delta<\delta_2$,
where $\delta_1<1<\delta_2$ are the two roots of $d(\delta)=J(u)$.
\end{lemma}
{\bf Proof}
If the sign of $I_\delta(u)$ changed for $\delta_1<\delta<\delta_2$,
then there exists $\delta_0\in(\delta_1,\delta_2)$ such that $I_{\delta_0}(u)=0$.
Thus $u\in\mathscr{N}_{\delta_0}$ and $d(\delta_0)\leq J(u)$.
According to Lemma \ref{lem3}, $d(\delta_0)>d(\delta_1)=d(\delta_2)=J(u)$,
which is a contradiction.
$\hfill\Box$

\begin{lemma}\label{lem5}
Assume that $u$ is a weak solution of \eqref{1.1} with $0<J(u_0)<d$
on $\Omega\times[0,T)$, $\delta_1<1<\delta_2$ are two roots of $d(\delta)=J(u_0)$.
\\
{\rm(i)} If $I(u_0)>0$, then $u(x,t)\in W_\delta$, $\delta_1<\delta<\delta_2$, $0<t<T$.
\\
{\rm(ii)} If $I(u_0)<0$, then $u(x,t)\in V_\delta$, $\delta_1<\delta<\delta_2$, $0<t<T$.
\end{lemma}
{\bf Proof}
(i) We first prove $u_0(x)\in W_\delta$ with $\delta_1<\delta<\delta_2$.
On the one hand, since $I(u_0)>0$ and Lemma \ref{lem4}, we have $I_{\delta}(u_0)>0$.
On the other hand, Lemma \ref{lem3} leads to $J(u_0)=d(\delta_1)=d(\delta_2)<d(\delta)$
with $\delta_1<\delta<\delta_2$. In what follows we prove that $u(x,t)\in W_\delta$
with $\delta_1<\delta<\delta_2$ on $0<t<T$. Suppose that there are $t_0\in (0,T)$
and $\delta_0\in (\delta_1,\delta_2)$ such that $u\in W_\delta$, $\delta_1<\delta<\delta_2$,
$0<t<t_0$, $u(x,t_0)\in \partial W_{\delta_0}$, then we can get
\begin{center}
$I_{\delta_0}(u(t_0))=0$, $\|\nabla u\|_p\neq 0$ or $J(u(t_0))=d(\delta_0)$.
\end{center}
Due to $\frac{d}{dt}J(u)\leq0$, then $J(u(t_0))\leq J(u_0)<d(\delta_0)$.
We only need to consider the first case, namely $u(t_0)\in \mathscr{N}_{\delta_0}$,
which indicates $J(u(t_0))\geq d(\delta_0)$. This is a contradiction.

(ii) The proof is similar to (i).
$\hfill\Box$

\section{$J(u_0)<d$}\label{sec1}

\qquad In this section, we establish the global existence and the blow-up properties
of the weak solution to \eqref{1.1} under the condition $J(u_0)<d$.

\begin{theorem}\label{th3.1}
Let $u_0\in W_{0}^{1,p}(\Omega)$ with $J(u_0)<d$ and $I(u_0)>0$.
Then \eqref{1.1} admits a unique global weak solution satisfying
$u\in L^\infty(0,\infty; W_{0}^{1,p}(\Omega))$ with $u_t\in L^2(0,\infty;L^2(\Omega))$.
Further, there exists $C>0$ such that $\|u\|_2^2+k\|\nabla u\|_2^2
\leq\left[(\|u_0\|_2^2+k\|\nabla u_0\|_2^2)^{1-p}+Ct\right]^{-\frac{1}{p-1}}$.
\end{theorem}
{\bf Proof} According to \eqref{3}, $I(u_0)>0$ and $2p<q^-+1$, we have $J(u_0)>0$.

{\bf Step1:} Global existence.

Let $\{\phi_j(x)\}_{j=1}^\infty$ be the orthogonal base in $W_{0}^{1,p}(\Omega)$,
which is also orthogonal in $L^{2}(\Omega)$.
Construct the approximate solution $u^m(x,t)$ of \eqref{1.1} as follows
\begin{align*}
u^m(x,t)=\sum_{j=1}^{m}\alpha_j^m(t)\phi_j(x),\qquad\alpha_j^m(t)=(u^m,\phi_j),\qquad m=1,2,...
\end{align*}
which satisfy
\begin{align}
\label{3.1.1}
&(u_{t}^m, \phi_j)+k(\nabla u_{t}^m, \nabla\phi_j)
+M(\|\nabla u^m\|^p_p)(|\nabla u^m|^{p-2}\nabla u^m, \nabla\phi_j)
=(|u^m|^{q-1}u^m\log|u^m|, \phi_j),
\\
\label{3.1.2}
&u^m(x,0)=\sum_{j=1}^{m}\alpha_j^m(0)\phi_j(x)\rightarrow u_0(x)\quad
\hbox{in}\quad W_{0}^{1,p}(\Omega).
\end{align}
Multiplying \eqref{3.1.1} by $\frac{d}{dt}\alpha_j^m(t)$, summing for $j$
from $1$ to $m$ and integrating with respect to time, we can obtain
\begin{align*}
J(u^m(x,0))=J(u^m(x,t))+\int_0^t(\|u_\tau^m\|_2^2+k\|\nabla u_\tau^m\|_2^2)d\tau,
\quad\forall t>0.
\end{align*}
By \eqref{3.1.2}, we have $J(u^m(x,0))\rightarrow J(u_0)<d$.
Hence for sufficiently large $m$, there holds
\begin{align*}
J(u^m(x,t))+\int_0^t(\|u_\tau^m\|_2^2+k\|\nabla u_\tau^m\|_2^2)d\tau=J(u^m(x,0))<d,
\quad\forall t>0.
\end{align*}
From \eqref{3.1.2} again, we have $I(u^m(x,0))\rightarrow I(u_0)>0$.
Hence for sufficiently large $m$, there holds $u^m(x,0)\in W$.
Then by Lemma \ref{lem5}, $u^m(x,t)\in W$ and
\begin{align*}
d>&\int_0^t(\|u_\tau^m\|_2^2+k\|\nabla u_\tau^m\|_2^2)d\tau
+\frac{a(q+1-p)}{p(q+1)}\|\nabla u^m\|^p_p
\\
&+\frac{b(q+1-2p)}{2p(q+1)}\|\nabla u^m\|^{2p}_p+\frac{1}{(q+1)^2}\|u^m\|_{q+1}^{q+1},
\quad\forall t>0.
\end{align*}
Thus
$$
\begin{aligned}
&\int_0^t(\|u_\tau^m\|_2^2+k\|\nabla u_\tau^m\|_2^2)d\tau<d,
\quad\|\nabla u^m\|^p_p<\frac{dp(q+1)}{a(q+1-p)},
\\
&\|u^m\|_{q+1}^{q+1}<d(q+1)^2.
\end{aligned}
$$
Then there exists a positive constant $C$ such that
$\|M(\|\nabla u^m\|_p^p)|\nabla u^m|^{p-2}\cdot\nabla u^m\|_{\frac{p}{p-1}}<C$.
Combining the Sobolev embedding inequality and
$\inf\{x^q\log x, x\in(0,1)\}=-\frac{1}{eq}$ with $q>0$, we have
\begin{align*}
\int_\Omega{\left|\left(|u^m|^{q}\log|u^m|\right)\right|^{\frac{q+2}{q+1}}}dx
=&\int_{\{x\in\Omega;|u^m|\leq1\}}{\left|\left(|u^m|^{q}\log|u^m|\right)\right|^{\frac{q+2}{q+1}}}dx
\\
&+\int_{\{x\in\Omega;|u^m|>1\}}{\left|\left(|u^m|^{q}\log|u^m|\right)\right|^{\frac{q+2}{q+1}}}dx
\\
\leq&\left(\frac{1}{eq}\right)^{\frac{q+2}{q+1}}\cdot|\Omega|+\|u^m\|_{q+2}^{q+2}
\\
\leq&\left(\frac{1}{eq}\right)^{\frac{q+2}{q+1}}\cdot|\Omega|+S^{q+2}\cdot
\left(\frac{pd(q+1)}{a(q+1-p)}\right)^{\frac{q+2}{p}}, \forall t>0.
\end{align*}
By the diagonal method and Aubin-Lion's compactness embedding theorem,
there exist $u$ and a subsequence of $\{u^m\}_{m=1}^{\infty}$
(still represented by $\{u^m\}_{m=1}^{\infty}$) such that
\begin{center}
$u_t^m\rightharpoonup u_t$ in $L^2(0,\infty;L^2(\Omega))$,
\\
$u^m\mathop{\rightharpoonup}\limits^* u$
in $L^\infty(0,\infty;W_0^{1,p}(\Omega))$,
\\
$u^m\rightarrow u$ strongly in $L^2(\Omega\times(0,T))$,
a.e. in $\Omega\times(0,T)$,
\\
$|u^m|^{q-1}u^m\cdot \log|u^m|\mathop{\rightharpoonup}\limits^* |u|^{q-1}u\cdot \log|u|$
in $L^\infty(0,\infty;L^{\frac{q+2}{q+1}}(\Omega))$,
\\
$M(\|\nabla u^m\|_p^p)|\nabla u^m|^{p-2}\cdot\nabla u^m\mathop{\rightharpoonup}\limits^* \xi$
in $L^\infty(0,\infty;L^{\frac{p}{p-1}}(\Omega))$.
\end{center}
Similar to the process of \cite{CaoEJDE, HanMMA}, we can prove
$\xi=M(\|\nabla u\|_p^p)|\nabla u|^{p-2}\nabla u$.
For fixed $j$, let $m\rightarrow +\infty$ in \eqref{3.1.1} to get
\begin{center}
$(u_t, \phi_j)+k(\nabla u_t, \nabla\phi_j)
+M(\|\nabla u\|^p_p)(|\nabla u|^{p-2}\nabla u, \nabla\phi_j)
=(|u|^{q-1}u\log|u|, \phi_j)$.
\end{center}
Then from Definition \ref{def2.1}, $u(x,t)$ is a global weak solution of \eqref{1.1}.

{\bf Step2:} Uniqueness.

Assume \eqref{1.1} has two global weak solution $u$ and $v$.
Set $w=u-v$, then $w$ satisfies
$$
\begin{aligned}
&\frac{1}{2}\frac{d}{dt}\int_{\Omega}w^2dx
+\frac{k}{2}\frac{d}{dt}\int_{\Omega}|\nabla w|^2dx
+M(\|\nabla u\|_p^p)\|\nabla u\|_p^p+M(\|\nabla v\|_p^p)\|\nabla v\|_p^p
\\
=&M(\|\nabla u\|_p^p)\int_{\Omega}|\nabla u|^{p-2}\nabla u\nabla vdx
+M(\|\nabla v\|^p)\int_{\Omega}|\nabla v|^{p-2}\nabla v\nabla udx
\\
&+\int_{\Omega}\left(q\left|\theta u
+(1-\theta)v\right|^{q-1}\log\left|\theta u+(1-\theta)v\right|
+\left|\theta u+(1-\theta)v\right|^{q-2}\cdot(\theta u+(1-\theta)v)\right)w^2dx
\end{aligned}
$$
with $\theta\in(0,1)$ and $w(x,0)=0$. Using the Young inequality, we can get
$$
\begin{aligned}
&\frac{1}{2}\frac{d}{dt}\int_{\Omega}w^2dx
+\frac{k}{2}\frac{d}{dt}\int_{\Omega}|\nabla w|^2dx
\\
\leq&\frac{1}{2}\frac{d}{dt}\int_{\Omega}w^2dx
+\frac{k}{2}\frac{d}{dt}\int_{\Omega}|\nabla w|^2dx
+\frac{1}{p}\left(M(\|\nabla u\|^p)-M(\|\nabla v\|^p)\right)\left(\|\nabla u\|^p-\|\nabla v\|^p\right)
\\
\leq&\int_{\Omega}\left(q\left|\theta u+(1-\theta)v\right|^{q-1}\log\left|\theta u+(1-\theta)v\right|
+\left|\theta u+(1-\theta)v\right|^{q-2}\cdot(\theta u+(1-\theta)v)\right)w^2dx.
\end{aligned}
$$
Therefore, by the Gronwall inequality, the boundedness of $u$ and $v$ and
$\inf\{x^q\log x, x\in(0,1)\}=-\frac{1}{eq}$ with $q>0$, we have $u=v$.

{\bf Step3:} Progressive estimation.

According to $u_0\in W$ and Lemma \ref{lem5}, we have
$u(x,t)\in W_\delta$, $\delta_1<\delta<\delta_2$,
where $\delta_1<1<\delta_2$ are two roots of $d(\delta)=J(u_0)$.
Furthermore, from the H\"{o}lder inequality and the Poincar\'{e} inequality,
there exist positive constants $C^*$ and $C_*$ such that
\begin{align*}
\frac{1}{2}\frac{d}{dt}(\|u\|_2^2+k\|\nabla u\|_2^2)
&=-I_{\delta_1}(u)+a(\delta_1-1)\|\nabla u\|_p^p+b(\delta_1-1)\|\nabla u\|_p^{2p}
\\
&\leq b(\delta_1-1)\|\nabla u\|_p^{2p}
\\
&\leq\frac{b(\delta_1-1)}{C^{*2p}}\|\nabla u\|_2^{2p}
\\
&\leq(\delta_1-1)\gamma(\|u\|_2^{2p}+k^p\|\nabla u\|_2^{2p})
\end{align*}
with $\gamma=\min\left\{\frac{b}{2k^pC^{*2p}}, \frac{bC_*^{2p}}{2C^{*2p}}\right\}$.
Since there exists $K_p>0$ for each $p$, such that $K_p(a^p+b^p)\geq(a+b)^p$
with non-negative $a$ and $b$, then
\begin{equation}\label{cao1}
\frac{1}{2}\frac{d}{dt}(\|u\|_2^2+k\|\nabla u\|_2^2)
\leq(\delta_1-1)\frac{\gamma}{K_p}\left(\|u\|_2^2+k\|\nabla u\|_2^2\right)^p,
\end{equation}
which implies $\|u\|_2^2+k\|\nabla u\|_2^2
\leq\left[\left(\|u_0\|_2^2+k\|\nabla u_0\|_2^2\right)^{1-p}
+(1-\delta_1)(p-1)\frac{2\gamma}{K_p}t\right]^{-\frac{1}{p-1}}$.
$\hfill\Box$

\begin{theorem}\label{th3.2}
Let $u_0\in W_0^{1,p}(\Omega)$ with $J(u_0)<d$ and $I(u_0)<0$.
Then the weak solution of \eqref{1.1} blows up in finite time,
namely there exists $T>0$, such that
$$
\lim\limits_{{t\rightarrow T^-}}\int_0^t(\|u\|_2^2+k\|\nabla u\|_2^2)d\tau=+\infty.
$$
\end{theorem}
{\bf Proof} Assume $u$ is a global solution of \eqref{1.1}. Let
\begin{align*}
H(t)=\int_0^t(\|u\|_2^2+k\|\nabla u\|_2^2)d\tau
+(T^*-t)(\|u_0\|_2^2+k\|\nabla u_0\|_2^2),\quad t\in [0,T^*],
\end{align*}
where $T^*$ is a sufficiently large time.
Then $H(t)\geq 0$ with $t\in [0,T^*]$. By a direct computation, we can get
\begin{align}\label{cao22}
&H'(t)=\|u\|_2^2+k\|\nabla u\|_2^2-\|u_0\|_2^2-k\|\nabla u_0\|_2^2,
\\
\label{cao23}
&H''(t)=2(u_t,u)+2k(\nabla u_t,\nabla u)=-2I(u),
\end{align}
and
\begin{align*}
(H'(t))^2&=4\left[\int_0^t((u_\tau,u)+k(\nabla u_\tau,\nabla u))d\tau\right]^2
\\
&\leq 4\left[\int_0^t(\|u_\tau\|_2^2+k\|\nabla u_\tau\|_2^2)d\tau\right]
\left[\int_0^t(\|u\|_2^2+k\|\nabla u\|_2^2)d\tau\right]
\\
&\leq 4H(t)\left[\int_0^t(\|u_\tau\|_2^2+k\|\nabla u_\tau\|_2^2)d\tau\right].
\end{align*}
Therefore we can deduce that
\begin{align}\label{cao25}
H''(t)H(t)-\frac{q+1}{2}(H'(t))^2
&\geq H(t)\left[-2I(u)-2(q+1)\int_0^t{\|u_\tau\|_2^2+k\|\nabla u_\tau\|_2^2}d\tau\right].
\end{align}
Set
$$
\xi(t)=-2I(u)-2(q+1)\int_0^t(\|u_\tau\|_2^2+k\|\nabla u_\tau\|_2^2)d\tau,
$$
which with \eqref{3} and the definition of $J(u)$ indicates
\begin{align*}
\xi(t)
=-2(q+1)J(u_0)+\frac{2a(q+1-p)}{p}\|\nabla u\|_p^p
+\frac{b(q+1-2p)}{p}\|\nabla u\|_p^{2p}
+\frac{2}{q+1}\|u\|_{q+1}^{q+1}.
\end{align*}
When $J(u_0)\leq 0$, then \eqref{cao2} means $J(u)\leq0$,
which with \eqref{3} leads to $I(u)<0$.
By Lemma \ref{lem2}, there holds $\|\nabla u\|_p>r(1)$. Therefore
\begin{equation}
\label{cao24}
\xi(t)>\sigma_1>0\quad \hbox{with}\quad \sigma_1=\frac{2a(q+1-p)}{p}r^p(1).
\end{equation}
When $0<J(u_0)<d$ and $I(u_0)<0$, then Lemma \ref{lem5} implies
$I_{\delta_2}(u)\leq 0$ and $\|\nabla u\|_p\geq r(\delta_2)>0$
with $\delta_1<1<\delta_2$ being the two roots of $J(u_0)=d(\delta)$.
Thus from \eqref{cao23}, we find that
\begin{align*}
H''(t)&=2a(\delta_2-1)\|\nabla u\|_p^p+2b(\delta_2-1)\|\nabla u\|_p^{2p}-2I_{\delta_2}(u)
\\
&\geq 2a(\delta_2-1)r^p(\delta_2),
\end{align*}
which with \eqref{cao22} guarantees
$$
\|u\|_2^2+k\|\nabla u\|_2^2\geq H'(t)\geq 2a(\delta_2-1)r^p(\delta_2)t.
$$
Thus there exists a $T_*>0$ such that \eqref{cao24} is established for $t\geq T_*$.
\\
Substituting \eqref{cao24} into \eqref{cao25}, we can deduce that
\begin{align*}
H''(t)H(t)-\frac{q+1}{2}(H'(t))^2>\sigma_1H(t).
\end{align*}
Then
\begin{align*}
[H^{\frac{1-q}{2}}(t)]''
\leq\frac{\sigma_1(1-q)}{2}[H^{\frac{1-q}{2}}(t)]^\frac{q+1}{q-1},
\quad t\in[T_*,T^*].
\end{align*}
Let $y(t)=H^{\frac{1-q}{2}}(t)$,
\begin{align*}
y''(t)\leq\frac{\sigma_1(1-q)}{2}[y(t)]^\frac{q+1}{q-1},
\quad t\in[T_*,T^*].
\end{align*}
Then there is $T\in (T_*,T^*)$ such that
$\lim\limits_{{t\rightarrow T^-}}y(t)=0$,
which means $\lim\limits_{{t\rightarrow T^-}}H(t)=+\infty$.
$\hfill\Box$


\section{$J(u_0)=d$}

\qquad In this section, we consider the behavior of solution to \eqref{1.1}
under the condition $J(u_0)=d$ and also give the threshold result
of global existence and finite time blowing-up.

\begin{lemma}\label{lem6}
If $u_0\in W_{0}^{1,p}(\Omega)$, $J(u_0)=d$, $I(u_0)>0$, then $W$ is an invariant set.
If $u_0\in W_{0}^{1,p}(\Omega)$, $J(u_0)=d$, $I(u_0)<0$, then $V$ is an invariant set.
\end{lemma}
{\bf Proof} Let $T$ be the maximum existence time of the solution.
If there exists $t_0\in(0,T)$ such that $I(u)>0$, $t\in[0,t_0)$ and $I(u(t_0))=0$.
Due to $-I(u)=(u_t,u)+k(\nabla u_t,\nabla u)<0$,
we get $\int_0^t(\|u_\tau\|_2^2+k\|\nabla u_\tau\|_2^2)d\tau>0$, $t\in(0,t_0)$.
Then
\begin{align*}
J(u(t_0))=J(u_0)-\int_0^{t_0}{\|u_\tau\|_2^2+k\|\nabla u_\tau\|_2^2}d\tau<d.
\end{align*}
It is known from $I(u(t_0))=0$ and the definition of $d$ \eqref{cao12}
that $J(u(t_0))\geq d$, which is a contradiction. Using the same method,
we can prove the second part of this lemma.
$\hfill\Box$

\begin{theorem}\label{th4.1}
Let $u_0\in W_{0}^{1,p}(\Omega)$, $J(u_0)=d$ and $I(u_0)\geq 0$.
Then \eqref{1.1} admits a unique global weak solution satisfies
$u\in L^\infty(0,\infty;W_{0}^{1,p}(\Omega))$,
$u_t\in L^2(0,\infty;L^2(\Omega))$ and $I(u)\geq 0$.
Moreover if $I(u)>0$, then there exists a constant $C>0$ and $t_0>0$ such that
$\|u\|_2^2+k\|\nabla u\|_2^2
\leq\left[(\|u(t_0)\|_2^2+k\|\nabla u(t_0)\|_2^2)^{1-p}+C(t-t_0)\right]^{-\frac{1}{p-1}}$.
Otherwise the solution vanishes in a limited time.
\end{theorem}
{\bf Proof} Since $J(u_0)=d$, $u_0>0$. Set $\lambda_s=1-\frac{1}{s}$, $s=1,2,...$
and consider the following initial value problem
\begin{eqnarray*}
\left\{
    \begin{array}{llll}
        \displaystyle u_{t}-k\Delta u_{t}-M(\|\nabla u\|_p^p)\Delta_p u=|u|^{q-1}u\log|u|,
        \quad&&(x,t)\in\Omega\times(0,T),
        \\
        \displaystyle u(x,t)=0, &&(x,t)\in\partial\Omega\times(0,T),
        \\
        \displaystyle u(x,0)=\lambda_s u_0(x),      && x\in\Omega.
    \end{array}
\right.
\end{eqnarray*}
According to $I(u_0)\geq 0$ and Lemma \ref{lem1}, there exists a unique
$\lambda^*\geq 1$ such that $I(\lambda^*u_0)=0$. Notice that
$\lambda_s<1\leq\lambda^*$, then $I(\lambda_su_0)>0$, $J(\lambda_su_0)<J(u_0)=d$.
By Theorem \ref{th3.1} and Lemma \ref{lem5}, for any $s$,
there exists a unique global weak solution
$u^s\in L^\infty(0,\infty;W_{0}^{1,p}(\Omega))$,
$u_t^s\in L^2(0,\infty;L^2(\Omega))$ such that $u^s\in W$ and
$$
\int_0^t(\|u_\tau^s\|_2^2+k\|\nabla u_\tau^s\|_2^2)d\tau+J(u^s)=J(\lambda_su_0)<d,
\quad 0\leq t<+\infty.
$$
Since $I(u^s)>0$, then we have
\begin{align*}
\int_0^t(\|u_\tau^s\|_2^2+k\|\nabla u_\tau^s\|_2^2)d\tau
+\frac{a(q+1-p)}{p(q+1)}\|\nabla u^s\|_p^p
+\frac{b(q+1-2p)}{2p(q+1)}\|\nabla u^s\|_p^{2p}
+\frac{1}{(q+1)^2}\|u^s\|_{q+1}^{q+1}<d.
\end{align*}
Therefore
\begin{align*}
&\int_0^t(\|u_\tau^s\|_2^2+k\|\nabla u_\tau^s\|_2^2)d\tau<d,
\\
&\|\nabla u^s\|_p^p<\frac{dp(q+1)}{a(q+1-p)},
\\
&\|u^s\|_{q+1}^{q+1}<d(q+1)^2.
\end{align*}
Similar to the proof of Theorem \ref{th3.1}, \eqref{1.1} has
a unique global weak solution $u\in L^\infty(0,\infty;W_{0}^{1,p}(\Omega))$,
$u_t\in L^2(0,\infty;L^2(\Omega))$ with $I(u)\geq 0$ and $J(u)\leq d$.
\\
If $I(u)>0$, $0<t<+\infty$, it can be seen from
$\frac{d}{dt}(\|u\|_2^2+k\|\nabla u\|_2^2)=-2I(u)<0$ that $u_t\neq0$.
Further there exists $t_0>0$ such that
$$
0<J(u(t_0))=J(u_0)-\int_0^{t_0}(\|u_\tau\|_2^2+k\|\nabla u_\tau\|_2^2)d\tau=d_1<d.
$$
Taking $t_0$ as the initial time, as it can be seen from Lemma \ref{lem5},
$u\in W_\delta$, $\delta_1<\delta<\delta_2$ for $t>t_0$,
where $\delta_1$ and $\delta_2$ are the two roots of $d(\delta)=J(u(t_0))$.
Therefore we have $I_{\delta_1}(u)\geq0$ for $t>t_0$.
Similar to the proof of \eqref{cao1}, we have
\begin{align*}
\frac{1}{2}\frac{d}{dt}(\|u\|_2^2+k\|\nabla u\|_2^2)
\leq\frac{(\delta_1-1)\gamma}{K_p}(\|u\|_2^2+k\|\nabla u\|_2^2)^p
\end{align*}
with $\gamma=\min\left\{\frac{b}{2k^pC^{*2p}}, \frac{bC_*^{2p}}{2C^{*2p}}\right\}$.
Then we can get
\begin{align*}
   \|u\|_2^2+k\|\nabla u\|_2^2
   \leq\left[(\|u(t_0)\|_2^2+k\|\nabla u(t_0)\|_2^2)^{1-p}
   +\frac{2(1-\delta_1)(p-1)\gamma}{K_p}(t-t_0)\right]^{-\frac{1}{p-1}}.
\end{align*}
If there exists $t^*>0$ such that $I(u)>0$, $0<t<t^*$, $I(u(t^*))=0$, then
\begin{align*}
&\int_0^{t^*}(\|u_\tau\|_2^2+k\|\nabla u_\tau\|_2^2)d\tau>0,
\\
&J(u(t^*))=d-\int_0^{t^*}(\|u_\tau\|_2^2+k\|\nabla u_\tau\|_2^2)d\tau<d.
\end{align*}
By the definition of $d$, we know $u(t^*)=0$.
Thus for all $t>t^*$, we have $u=0$, which means the weak solution of \eqref{1.1}
distinguishes at a finite time.
$\hfill\Box$

\begin{theorem}\label{th4.2}
Let $u_0\in W_{0}^{1,p}(\Omega)$, $J(u_0)=d$ and $I(u_0)<0$.
Then the weak solution of \eqref{1.1} blows up in finite time,
namely there exists $T>0$ such that
$$
\lim\limits_{{t\rightarrow T^-}}\int_0^t(\|u\|_2^2+k\|\nabla u\|_2^2)d\tau=+\infty.
$$
\end{theorem}
{\bf Proof} Similar to Theorem \ref{th3.2}, we get
\begin{align*}
&H''(t)H(t)-\frac{q+1}{2}(H'(t))^2
\\
\geq&\left[\frac{2a(q+1-p)}{p}\|\nabla u\|_p^p+\frac{b(q+1-2p)}{p}\|\nabla u\|_p^{2p}
+\frac{2}{q+1}\|u\|_{q+1}^{q+1}-2(q+1)J(u_0)\right]H(t).
\end{align*}
From $J(u_0)=d>0$, $I(u_0)<0$ and Lemma \ref{lem6}, there exists $t_0>0$ such that $I(u(t))<0$, $0<t<t_0$.
Then \eqref{cao23} leads to $H''(t)>0$ and $\|u_t\|_2^2+k\|\nabla u_t\|_2^2\neq 0$ for $0<t<t_0$.
Therefore $J(u(t_0))=d-\int_0^{t_0}(\|u_\tau\|^2+k\|\nabla u_\tau\|^2)d\tau=d_1<d$.
We choose $t_0$ as the initial time and complete the proof according to Theorem \ref{th3.2}.
$\hfill\Box$

\section{Decay rate and Life span}

\qquad In this section, we prove the exponential decay of global solution
and establish life span estimation of finite time blow-up solution
without additional restriction on the initial data.

\begin{theorem}\label{th1.3}
Let $u_0(x)\in W_0^{1,p}(\Omega)$ with $J(u_0)<d$ and $I(u_0)>0$.
Then the global weak solutions of \eqref{1.1} satisfy
\begin{align*}
J(u)&<G(0)e^{-\frac{\alpha}{\tilde{\alpha}t}},
\\
\|\nabla u\|_p^p&<\frac{p(q+1)G(0)}{a(q+1-p)}e^{-\frac{\alpha}{\tilde{\alpha}}t},
\\
\|u\|_{q+1}^{q+1}
&<\left(\frac{S_1^pp(q+1)G(0)}{a(q+1-p)}\right)^{\frac{q+1}{p}}e^{-\frac{\alpha(q+1)}{\tilde{\alpha}p}t},
\end{align*}
where $G(0)=J(u_0)+\|u_0\|_2^2+k\|\nabla u_0\|_2^2$,
$\beta=\left(\frac{p(q+1)J(u_0)}{a(q+1-p)}\right)^\frac{q+1-p}{p}$,
$\tilde{\alpha}$ is a positive constant,
$\alpha=\frac{4ap(q+1)^2(1-\delta_3)}{a(q+1)(3q+3-2p-2p\delta_3)+2pS_1^{q+1}\beta}>0$,
$S_1$ is the Sobolev embedding constant with $\|u\|_{q+1}<S_1\|\nabla u\|_p$
and $\delta_1<\delta_3=\frac{1+\delta_1}{2}<1$ with $\delta_1$ in Lemma \ref{lem5}.
\end{theorem}
{\bf Proof} From Lemma \ref{lem5}, we get that $u\in W_\delta$,
where $\delta_1<\delta<\delta_2$ with $\delta_1<1<\delta_2$ satisfies $J(u_0)=d(\delta)$.
Take $\delta_1<\delta_3=\frac{1+\delta_1}{2}<1$, then $I_{\delta_3}(u)>0$.
Thus by $I(u)=I_{\delta_3}(u)+a(1-\delta_3)\|\nabla u\|_p^p
+b(1-\delta_3)\|\nabla u\|_p^{2p}$ and the Sobolev embedding inequality, we can get
\begin{align}
&a\|\nabla u\|_p^p<\frac{1}{1-\delta_3}I(u),\quad
b\|\nabla u\|_p^{2p}<\frac{1}{1-\delta_3}I(u),\nonumber
\\
\label{cao17}
&\|u\|_{q+1}^{q+1}<S_1^{q+1}\|\nabla u\|_p^{q+1}
<\frac{S_1^{q+1}\beta}{a(1-\delta_3)}I(u).
\end{align}
Invoke by the ideas in \cite{XuNA}, we set
\begin{equation}\label{cao18}
G(t)=J(u(t))+\|u\|_2^2+k\|\nabla u\|_2^2.
\end{equation}
Combining \eqref{cao17}, \eqref{cao18}, the H\"{o}lder inequality,
the Sobolev inequality and the assumption on the exponent $p$ and $q$,
there exists a positive constant $\tilde{\alpha}$ such that
\begin{equation}\label{cao19}
J(u(t))<G(t)<\tilde{\alpha}J(u(t)).
\end{equation}
By direct computation and \eqref{cao19}, we get
\begin{align*}
G'(t)=&-\|u_t\|_2^2-k\|\nabla u_t\|_2^2-2I(u)-\alpha J(u)
+\frac{\alpha}{q+1}I(u)+\frac{\alpha a(q+1-p)}{p(q+1)}\|\nabla u\|_p^p
\\
&+\frac{\alpha b(q+1-2p)}{2p(q+1)}\|\nabla u\|_p^{2p}
+\frac{\alpha}{(q+1)^2}\|u\|_{q+1}^{q+1}
\\
\leq&-\alpha J(u)+\left[\frac{\alpha}{q+1}-2
+\frac{\alpha(q+1-p)}{p(q+1)(1-\delta_3)}
+\frac{\alpha(q+1-2p)}{2p(q+1)(1-\delta_3)}
+\frac{\alpha S_1^{q+1}\beta}{a(1-\delta_3)(q+1)^2}\right]I(u).
\end{align*}
Choosing $\alpha=\frac{4ap(q+1)^2(1-\delta_3)}{a(q+1)(3q+3-2p-2p\delta_3)+2pS_1^{q+1}\beta}>0$,
which makes the coefficient of $I(u)$ being $0$, then from \eqref{cao19} we have
\begin{align*}
G'(t)\leq-\alpha J(u)\leq-\frac{\alpha}{\tilde{\alpha}}G(t).
\end{align*}
Thus
\begin{align*}
G(t)\leq G(0)e^{-\frac{\alpha}{\tilde{\alpha}}t},
\end{align*}
which with \eqref{cao19} leads to
\begin{align*}
J(u)<G(t)\leq G(0)e^{-\frac{\alpha}{\tilde{\alpha}}t}.
\end{align*}
Furthermore, using \eqref{3}, Lemma \ref{lem5} and the above inequality imply that
\begin{align*}
\|\nabla u\|_p^p&\leq\frac{p(q+1)}{a(q+1-p)}J(u)
<\frac{p(q+1)G(0)}{a(q+1-p)}e^{-\frac{\alpha}{\tilde{\alpha}}t},
\\
\|u\|_{q+1}^{q+1}
&<\left(\frac{S_1^pp(q+1)G(0)}{a(q+1-p)}\right)^{\frac{q+1}{p}}e^{-\frac{\alpha(q+1)}{\tilde{\alpha}p}t}.
\end{align*}
$\hfill\Box$

\begin{theorem}\label{th4.3}
Let $u_0\in W_{0}^{1,p}(\Omega)$, $J(u_0)\leq d$ and $I(u_0)<0$.
Then we have the following life span estimation of the blow-up
solution in Theorem \ref{th3.2} and Theorem \ref{th4.2}.
\\
{\rm(i)} If $J(u_0)<0$, then $T\leq\frac{\|u_0\|_2^2+k\|\nabla u_0\|_2^2}{(1-q^2)J(u_0)}$.
\\
{\rm(ii)} If $0\leq J(u_0)\leq d$, then
$T\leq\frac{4\|u_0\|_2^2+4k\|\nabla u_0\|_2^2}{(q-1)^2(\frac{a(q+1-p)}{p(q-1)}
\|\nabla u(t_0)\|_p^p+\frac{b(q+1-2p)}{2p(q-1)}\|\nabla u(t_0)\|_p^{2p}-J(u_0))}+t_0$,
where $t_0$ satisfies \eqref{cao3}.
\end{theorem}
{\bf Proof}
Let $\theta(t)=\frac{1}{2}\|u\|_2^2+\frac{k}{2}\|\nabla u\|_2^2$, $\eta(t)=-J(u)$.

(i) If $J(u_0)<0$, then we have $\theta(0)>0$, $\eta(0)>0$,
$\eta'(t)=\|u_t\|_2^2+k\|\nabla u_t\|_2^2\geq 0$ and $\eta(t)>0$.
From a simple computation, we can find that
\begin{align}
\theta'(t)&=-I(u)\nonumber
\\
&=-(q+1)J(u)+\frac{a(q+1-p)}{p}\|\nabla u\|_p^p
+\frac{b(q+1-2p)}{2p}\|\nabla u\|_p^{2p}+\frac{1}{q+1}\|u\|_{q+1}^{q+1}\nonumber
\\
\label{cao26}
&>(q+1)\eta(t)>0,
\\
\theta(t)\eta'(t)&=\frac{1}{2}(\|u\|_2^2
+k\|\nabla u\|_2^2)\cdot(\|u_t\|_2^2+k\|\nabla u_t\|_2^2)\nonumber
\\
&\geq\frac{1}{2}((u_t,u)+k(\nabla u_t,\nabla u))^2\nonumber
\\
&=\frac{1}{2}(\theta'(t))^2>\frac{q+1}{2}\theta'(t)\eta(t),\nonumber
\end{align}
which implies
\begin{align}\label{cao27}
\left[\eta(t)\cdot\theta^{\frac{-q-1}{2}}(t)\right]'
=\theta^{-1-\frac{q+1}{2}}(t)\left[\theta(t)\eta'(t)
-\frac{q+1}{2}\eta(t)\theta'(t)\right]>0.
\end{align}
Thus \eqref{cao26} and \eqref{cao27} lead to
\begin{align*}
0<\eta(0)\theta^{\frac{-q-1}{2}}(0)
\leq\eta(t)\theta^{\frac{-q-1}{2}}(t)
\leq\frac{1}{q+1}\theta'(t)\theta^{\frac{-q-1}{2}}(t)
=\frac{2}{1-q^2}\left[\theta^{\frac{1-q}{2}}(t)\right]',
\end{align*}
which further indicates that
\begin{align*}
0\leq\theta^{\frac{1-q}{2}}(t)
\leq\frac{1-q^2}{2}\theta^{\frac{-q-1}{2}}(0)\eta(0)t+\theta^{\frac{1-q}{2}}(0).
\end{align*}
Thus we can deduce that
$T\leq \frac{\|u_0\|_2^2+k\|\nabla u_0\|_2^2}{(1-q^2)J(u_0)}$.

(ii) For $0\leq J(u_0)\leq d$, if there exists $t^*>0$ such that $J(u(t^*))<0$,
we can get the upper bound estimation of $T$ by (i).
Therefore we only need to consider the case $0\leq J(u)\leq d$
for any $t\in (0,T)$. According to Lemma \ref{lem4} and Lemma \ref{lem6},
we have $I(u)<0$ with $0<t<T$, which means $\theta'(t)=-I(u)>0$.
Since $u$ blows up in finite time and $\theta(t)$ increases
with respect to time, there exists $0<t_0<T$ such that
\begin{align}\label{cao3}
&2(q+1)J(u_0)\nonumber
\\
&<\min\limits_{t\in[t_0,T)}\left(\frac{2a(q+1-p)}{p}\|\nabla u(t)\|_p^p
+\frac{b(q+1-2p)}{p}\|\nabla u(t)\|_p^{2p}
+\frac{2}{q+1}\|u(t)\|_{q+1}^{q+1}\right).
\end{align}
Let
\begin{align*}
F(t)=\int_{t_0}^t(\|u\|_2^2+k\|\nabla u\|_2^2)d\tau
+(T-t)(\|u_0\|_2^2+k\|\nabla u_0\|_2^2)+\beta((t-t_0)+\sigma)^2
\end{align*}
with $t_0\leq t<T$ and positive constants $\beta$ and $\sigma$ to be determined later.
By a direct computation, we can get
\begin{align*}
F'(t)=&\int_{t_0}^t{\frac{d}{d\tau}(\|u\|_2^2+k\|\nabla u\|_2^2)}d\tau
+2\beta((t-t_0)+\sigma),
\\
F''(t)=&\frac{d}{dt}(\|u\|_2^2+k\|\nabla u\|_2^2)+2\beta
\\
=&-2(q+1)J(u(t_0))
+2(q+1)\int_{t_0}^t(\|u_\tau\|_2^2+k\|\nabla u_\tau\|_2^2)d\tau
+\frac{2a(q+1-p)}{p}\|\nabla u\|_p^p
\\
&+\frac{b(q+1-2p)}{p}\|\nabla u\|_p^{2p}
+\frac{2}{q+1}\|u\|_{q+1}^{q+1}+2\beta.
\end{align*}
Then $F(t_0)>0$, $F'(t_0)>0$ and $F'(t)>0$ for $t\in [t_0,T)$,
provided that $\beta$ or $\sigma$ is large enough.
For any $\rho>0$, we have
\begin{align*}
&F''(t)F(t)-\rho(F'(t))^2
\\
=&F''(t)F(t)+4\rho\Big[
\left(\int_{t_0}^t(\|u\|_2^2+k\|\nabla u\|_2^2)d\tau
+\beta((t-t_0)+\sigma)^2\right)
\left(\int_{t_0}^t(\|u_\tau\|_2^2+k\|\nabla u_\tau\|_2^2)d\tau
+\beta\right)
\\
&-\left(\int_{t_0}^t[(u_\tau,u)+k(\nabla u_\tau,\nabla u)]d\tau
+\beta((t-t_0)+\sigma)\right)^2
\\
&-(F(t)-(T-t)(\|u_0\|_2^2+k\|\nabla u_0\|_2^2))
\left(\int_{t_0}^t(\|u_\tau\|_2^2+k\|\nabla u_\tau\|_2^2)d\tau
+\beta\right)\Big]
\\
=&F''(t)F(t)+4\rho(T-t)(\|u_0\|_2^2+k\|\nabla u_0\|_2^2)
\left(\int_{t_0}^t(\|u_\tau\|_2^2+k\|\nabla u_\tau\|_2^2)d\tau
+\beta\right)
\\
&+4\rho\zeta(t)-4\rho F(t)\left(\int_{t_0}^t(\|u_\tau\|_2^2
+k\|\nabla u_\tau\|_2^2)d\tau+\beta\right),
\end{align*}
where
\begin{align*}
\zeta(t)=&\left(\int_{t_0}^t(\|u\|_2^2+k\|\nabla u\|_2^2)d\tau
+\beta((t-t_0)+\sigma)^2\right)
\left(\int_{t_0}^t(\|u_\tau\|_2^2+k\|\nabla u_\tau\|_2^2)d\tau
+\beta\right)
\\
&-\left(\int_{t_0}^t[(u_\tau, u)+k(\nabla u_\tau, \nabla u)]d\tau
+\beta((t-t_0)+\sigma)\right)^2\geq 0,
\quad t\in[t_0,T).
\end{align*}
Thus
\begin{align*}
&F''(t)F(t)-\rho(F'(t))^2
\\
\geq&F(t)\left[F''(t)-4\rho\beta
-4\rho\int_{t_0}^t(\|u_\tau\|_2^2+k\|\nabla u_\tau\|_2^2)d\tau\right]
\\
\geq&F(t)\left[-2(q+1)J(u(t_0))+2(q+1-2\rho)\int_{t_0}^t(\|u_\tau\|_2^2
+k\|\nabla u_\tau\|_2^2)d\tau\right.
\\
&\left.+\frac{2a(q+1-p)}{p}\|\nabla u\|_p^p
+\frac{b(q+1-2p)}{p}\|\nabla u\|_p^{2p}+\frac{2}{q+1}\|u\|_{q+1}^{q+1}
+2\beta-4\rho\beta\right].
\end{align*}
Take $\rho=\frac{q+1}{2}$ and using \eqref{cao2}, the above inequality can be reduced to
\begin{align*}
&F''(t)F(t)-\frac{q+1}{2}(F'(t))^2
\\
\geq&F(t)
\left[\frac{2a(q+1-p)}{p}\|\nabla u\|_p^p+\frac{b(q+1-2p)}{p}\|\nabla u\|_p^{2p}
+\frac{2}{q+1}\|u\|_{q+1}^{q+1}-2(q+1)J(u(t_0))-2(q+1)\beta\right].
\end{align*}
Since \eqref{cao3}, we can get
$$
F''(t)F(t)-\frac{q+1}{2}(F'(t))^2\geq 0\quad\hbox{with}\quad t_0<t<T,
$$
provided that
\begin{equation}\label{cao28}
\beta\in\left(0, \frac{a(q+1-p)}{p(q+1)}\|\nabla u(t_0)\|_p^p
+\frac{b(q+1-2p)}{2p(q+1)}\|\nabla u(t_0)\|_p^{2p}
+\frac{1}{(q+1)^2}\|u(t_0)\|_{q+1}^{q+1}-J(u(t_0))\right].
\end{equation}
Set $G(t)=F^{1-\frac{q+1}{2}}(t)$ with $t\in[t_0,T)$, then
\begin{align*}
G'(t)&=\left(1-\frac{q+1}{2}\right)F^{-\frac{q+1}{2}}(t)F'(t)\leq 0,
\\
G''(t)&=\left(1-\frac{q+1}{2}\right)\cdot
F^{-1-\frac{q+1}{2}}(t)\left[-\frac{q+1}{2}(F'(t))^2
+F''(t)F(t)\right]\leq 0,
\\
G(t)&\leq G(t_0)+G'(t_0)(t-t_0).
\end{align*}
Because of $G(t_0)>0$ and $G'(t_0)<0$, we have
\begin{align*}
t-t_0\leq-\frac{G(t_0)}{G'(t_0)}
=\frac{(T-t_0)(\|u_0\|_2^2+k\|\nabla u_0\|_2^2)+\beta\sigma^2}{(q-1)\beta\sigma}
\end{align*}
with $t\in(t_0,T)$.
For fixed $\beta_0$ satisfies \eqref{cao28}, taking
$\sigma\in\left(\frac{\|u_0\|_2^2+k\|\nabla u_0\|_2^2}{(q-1)\beta_0}, +\infty\right)$, then
\begin{align*}
T-t_0\leq\frac{\beta_0\sigma^2}{(q-1)\beta_0\sigma-(\|u_0\|_2^2+k\|\nabla u_0\|_2^2)}.
\end{align*}
Define
$$
T_{\beta_0}(\sigma)
=\frac{\beta_0\sigma^2}
{(q-1)\beta_0\sigma-(\|u_0\|_2^2+k\|\nabla u_0\|_2^2)}
$$
with $\sigma\in\left(\frac{\|u_0\|_2^2+k\|\nabla u_0\|_2^2}{(q-1)\beta_0}, +\infty\right)$.
We can find that $T_{\beta_0}(\sigma)$ takes the minimum at
$$
\sigma=\frac{2\|u_0\|_2^2+2k\|\nabla u_0\|_2^2}{(q-1)\beta_0},
$$
which indicates that
\begin{align*}
T-t_0
\leq\inf\limits_{\sigma\in\left(\frac{\|u_0\|_2^2+k\|\nabla u_0\|_2^2}{(q-1)\beta_0}, +\infty\right)}
T_{\beta_0}(\sigma)
=\frac{4\|u_0\|_2^2+4k\|\nabla u_0\|_2^2}{(q-1)^2\beta_0}.
\end{align*}
Combining the above inequality with \eqref{cao28}, we can get
\begin{align*}
T-t_0&\leq\inf\limits_{\beta_0\in\left(0,
\frac{a(q+1-p)}{p(q+1)}\|\nabla u(t_0)\|_p^p
+\frac{b(q+1-2p)}{2p(q+1)}\|\nabla u(t_0)\|_p^{2p}
+\frac{1}{(q+1)^2}\|u(t_0)\|_{q+1}^{q+1}-J(u_0)\right]}
\frac{4\|u_0\|_2^2+4k\|\nabla u_0\|_2^2}{(q-1)^2\beta_0}
\\
&=\frac{4\|u_0\|_2^2+4k\|\nabla u_0\|_2^2}
{(q-1)^2(\frac{a(q+1-p)}{p(q+1)}\|\nabla u(t_0)\|_p^p
+\frac{b(q+1-2p)}{2p(q+1)}\|\nabla u(t_0)\|_p^{2p}+\frac{1}{(q+1)^2}\|u(t_0)\|_{q+1}^{q+1}-J(u_0))}.
\end{align*}
$\hfill\Box$

\section{$J(u_0)>d$}

\qquad In this section, we investigate the conditions
that ensure the global existence or finite time blowing-up
of solution to \eqref{1.1}. Inspired by the ideas in \cite{XS,Zhoujun1,Han2,HanMMA},
we first introduce the following sets.
\begin{align*}
\mathcal{N}_+&=\{u\in W_{0}^{1,p}(\Omega): I(u)>0\},
\\
\mathcal{N}_-&=\{u\in W_{0}^{1,p}(\Omega): I(u)<0\},
\\
J^s&=\{u\in W_{0}^{1,p}(\Omega): J(u)<s\},\quad\hbox{for any}\,s>d,
\\
\mathcal{N}^s&=\mathcal{N}\cap J^s=\{u\in\mathcal{N}: J(u)<s\},
\\
\lambda_s&=\inf\{\|u\|_2^2+k\|\nabla u\|_2^2: u\in\mathcal{N}^s\},
\\
\Lambda_s&=\sup\{\|u\|_2^2+k\|\nabla u\|_2^2: u\in\mathcal{N}^s\},
\\
\mathcal{B}&=\{u_0\in W_{0}^{1,p}(\Omega):\,\hbox{the solution of}\,\eqref{1.1}\,\hbox{blows up in finite time}\},
\\
\mathcal{G}&=\{u_0\in W_{0}^{1,p}(\Omega):\,\hbox{the solution of}\,\eqref{1.1}\,\hbox{is global in time}\},
\\
\mathcal{G}_0&=\{u_0\in W_{0}^{1,p}(\Omega):\,u(t)\rightarrow0\,\hbox{in}\,W_{0}^{1,p}(\Omega), t\rightarrow+\infty\}.
\end{align*}

\begin{theorem}\label{th5.1}
Assume $J(u_0)>d$.
\\
{\rm(i)} If $u_0\in\mathcal{N}_+$,
$\|u_0\|_2^2+k\|\nabla u_0\|_2^2\leq\lambda_{J(u_0)}$, then $u_0\in\mathcal{G}_0$.
\\
{\rm(ii)} If $u_0\in\mathcal{N}_-$,
$\|u_0\|_2^2+k\|\nabla u_0\|_2^2\geq\Lambda_{J(u_0)}$, then $u_0\in\mathcal{B}$.
\end{theorem}
{\bf Proof}
(i) If $u_0\in\mathcal{N}_+$ and $\|u_0\|_2^2+k\|\nabla u_0\|_2^2\leq\lambda_{J(u_0)}$,
then we assert that $u(t)\in\mathcal{N}_+$, $0\leq t<T(u_0)$ with $T(u_0)$
being the maximum existence time of the solution.
Otherwise there exists $t_0\in(0, T(u_0))$ such that $u(t)\in\mathcal{N}_+$,
$0\leq t<t_0$ and $u(t_0)\in\mathcal{N}$. Furthermore, \eqref{cao2} indicates that
$J(u(t_0))<J(u_0)$, which with the definition of $J^s$ leads to $u(t_0)\in J^{J(u_0)}$.
Thus $u(t_0)\in\mathcal{N}^{J(u_0)}$. According to the definition of $\lambda_{J(u_0)}$,
we can get
\begin{equation}\label{cao11}
\|u(t_0)\|_2^2+k\|\nabla u(t_0)\|_2^2\geq\lambda_{J(u_0)}.
\end{equation}
It can be seen from $u(t)\in \mathcal{N}_+$ with $0\leq t<t_0$ that
$$
I(u)=-\frac{1}{2}\frac{d}{dt}(\|u\|_2^2+k\|\nabla u\|_2^2)>0,\quad 0<t<t_0.
$$
Then we have
$$
\|u(t_0)\|_2^2+k\|\nabla u(t_0)\|_2^2<\|u_0\|_2^2
+k\|\nabla u_0\|_2^2\leq\lambda_{J(u_0)},
$$
which contradicts \eqref{cao11}.
Hence $u(t)\in\mathcal{N}_+$, $0\leq t<T(u_0)$.
\\
Using \eqref{3} and \eqref{cao2}, there holds
\begin{align*}
J(u_0)\geq J(u)
&={\frac{1}{q+1}}I(u)+({\frac{a}{p}}-{\frac{a}{q+1}})\|\nabla u\|_p^p
+({\frac{b}{2p}}-{\frac{b}{q+1}})\|\nabla u\|_p^{2p}
+\frac{1}{(q+1)^2}\|u\|_{q+1}^{q+1}
\\
&>\frac{a(q+1-p)}{p(q+1)}\|\nabla u\|_p^p,
\end{align*}
which means $\|\nabla u\|_p^p\leq\frac{p(q+1)J(u_0)}{a(q+1-p)}$
and further $T(u_0)=+\infty$.
Define the $\omega-\text{limit}$ set of $u_0$ by
$\omega(u_0)=\bigcap\limits_{t\geq 0}\overline{\{u(\cdot,s):s\geq t\}}$.
Then for any $\omega\in\omega(u_0)$, we have
$$
\|\omega\|_2^2+k\|\nabla \omega\|_2^2
<\|u_0\|_2^2+k\|\nabla u_0\|_2^2
\leq\lambda_{J(u_0)},
\quad J(\omega)\leq J(u_0).
$$
So that $\omega(u_0)\cap\mathcal{N}=\emptyset$,
which with the convergence result in \cite{CaoAKN} leads to $\omega(u_0)=\{0\}$,
namely $u_0\in\mathcal{G}_0$.

(ii) If $u_0\in\mathcal{N}_-$, $\|u_0\|_2^2+k\|\nabla u_0\|_2^2\geq\Lambda_{J(u_0)}$,
then similar to (i), we can get $u(t)\in\mathcal{N}_-$, $u(t)\in J^{J(u_0)}$
for $0\leq t<T(u_0)$. If $T(u_0)=\infty$, then for any $\omega\in\omega(u_0)$,
we conclude that
$$
\|\omega\|_2^2+k\|\nabla \omega\|_2^2>\Lambda_{J(u_0)}, J(\omega)\leq J(u_0).
$$
Then $\omega(u_0)\cap\mathcal{N}=\emptyset$,
which with the convergence result in \cite{CaoAKN} leads to $\omega(u_0)=\{0\}$.
However due to $u\in\mathcal{N}_-$, we have
$$
a\|\nabla u\|_p^p<a\|\nabla u\|_p^p+b\|\nabla u\|_p^{2p}<\int_\Omega{|u|^{q+1}log|u|}dx
\leq\|u\|_{q+2}^{q+2}<S^{q+2}\|\nabla u\|_p^{q+2},
$$
which means $\|\nabla u\|_p\geq\left(\frac{a}{S^{q+2}}\right)^{\frac{1}{q+2-p}}$.
It is a contradiction. Then $T(u_0)<+\infty$ and $u_0\in\mathcal{B}$.
$\hfill\Box$

\begin{theorem}\label{th1.6}
$\lambda_s\geq\left\{
\begin{aligned}
&\left[\frac{a}{\beta^{q+2}}\kappa^{p-\theta(q+2)}\right]^{\frac{2}{(1-\theta)(q+2)}}, &p>\frac{n}{n+2}(q+2),
\\
&\left[\frac{a}{\beta^{q+2}}\tilde{\kappa}^{\frac{p-\theta(q+2)}{p}}\right]^{\frac{2}{(1-\theta)(q+2)}}, &p<\frac{n}{n+2}(q+2),
\end{aligned}
\right.$
\\
and
$$
\Lambda_s\leq(1+k)|\Omega|^{\frac{p-2}{p}}\tilde{\kappa}^2,
$$
where $\theta$ satisfy $\theta(\frac{1}{2}-\frac{1}{p}+\frac{1}{n})=\frac{1}{2}-\frac{1}{q+2}$,
$\kappa$ is the unique positive solution of $f(y)=d$ and
$\tilde{\kappa}$ is the unique positive solution of $f(y)=s$ with
\begin{equation}\label{cao13}
f(y)=\frac{b(q+1-2p)}{2p(q+1)}y^{2p}+\frac{a(q+1-p)}{p(q+1)}y^p
+\frac{S_1^{q+1}}{(q+1)^2}y^{q+1},
\quad y\in\mathbb{R}.
\end{equation}

\end{theorem}
{\bf Proof}
By the definition of the potential well depth $d$ \eqref{cao12}, we have
\begin{align*}
d&=\inf\limits_{u\in\mathscr{N}}J(u)
\\
&=\inf\limits_{u\in\mathscr{N}}\left[\frac{a(q+1-p)}{p(q+1)}\|\nabla u\|_p^p+
\frac{b(q+1-2p)}{2p(q+1)}\|\nabla u\|_p^{2p}+\frac{1}{(q+1)^2}\|u\|_{q+1}^{q+1}\right]
\\
&\leq\inf\limits_{u\in\mathscr{N}}\left[\frac{a(q+1-p)}{p(q+1)}\|\nabla u\|_p^p+
\frac{b(q+1-2p)}{2p(q+1)}\|\nabla u\|_p^{2p}+\frac{S_1^{q+1}}{(q+1)^2}\|\nabla u\|_p^{q+1}\right]
\\
&=\inf\limits_{u\in\mathscr{N}}f(\|\nabla u\|_p),
\end{align*}
where $f(\cdot)$ is given in \eqref{cao13}.
Since $f(\cdot)$ is strictly increasing on $[0,+\infty)$ and $f(0)=0$,
there exists a unique $\kappa$ such that $f(\kappa)=d$.
Then for any $u\in\mathscr{N}$, there is
\begin{equation}\label{cao14}
\|\nabla u\|_p\geq\kappa>0.
\end{equation}

By the Gagliardo--Nirenberg inequality \cite{BrezisBOOK}, we get
\begin{align*}
\|u\|_{q+2}\leq\beta\|u\|_2^{(1-\theta)}\|\nabla u\|_p^\theta,
\end{align*}
where $\beta$ is a positive constant and
$\theta(\frac{1}{2}-\frac{1}{p}+\frac{1}{n})=\frac{1}{2}-\frac{1}{q+2}$.
Then it follows from the above inequality that for any $u\in\mathscr{N}$
\begin{align}
a\|\nabla u\|_p^p\leq\int_\Omega{|u|^{q+1}\log|u|}dx
&\leq\|u\|_{q+2}^{q+2}\leq\beta^{q+2}\|u\|_2^{(1-\theta)(q+2)}
\|\nabla u\|_p^{\theta(q+2)},\nonumber
\\
\label{cao15}
\text{which says}\quad a\|\nabla u\|_p^{p-\theta(q+2)}&\leq\beta^{q+2}\|u\|_2^{(1-\theta)(q+2)}.
\end{align}
Moreover by the definition of $\mathscr{N}^s$,
it is known that if $u\in\mathscr{N}^s$,
\begin{equation}\label{cao16}
\|\nabla u\|_p\leq \tilde{\kappa}=\left(\frac{ps(q+1)}{a(q+1-p)}\right)^{1/p}.
\end{equation}

For the lower bound of $\lambda_s$, we divide into two cases to discuss.
\\
Case 1: $p-\theta(q+2)>0$, namely $p>\frac{n}{n+2}(q+2)$,
then using \eqref{cao14} and \eqref{cao15}, we have
\begin{align*}
\lambda_s&=\inf\limits_{u\in\mathscr{N}^s}\left\{\|u\|_2^2+k\|\nabla u\|_2^2\right\}
\\
&\geq\inf\limits_{u\in\mathscr{N}}\left\{\|u\|_2^2+k\|\nabla u\|_2^2\right\}
\\
&\geq\inf\limits_{u\in\mathscr{N}
}\left[\frac{a}{\beta^{q+2}}\|\nabla u\|_p^{p-\theta(q+2)}\right]^{\frac{2}{(1-\theta)(q+2)}}
\\
&\geq\left[\frac{a}{\beta^{q+2}}\kappa^{p-\theta(q+2)}\right]^{\frac{2}{(1-\theta)(q+2)}}.
\end{align*}
\\
Case 2: $p-\theta(q+2)<0$, namely $p<\frac{n}{n+2}(q+2)$,
then using \eqref{cao15} and \eqref{cao16}, we have
\begin{align*}
\lambda_s&=\inf\limits_{u\in\mathscr{N}^s}\left\{\|u\|_2^2+k\|\nabla u\|_2^2\right\}
\\
&\geq\inf\limits_{u\in\mathscr{N}}
\left[\frac{a}{\beta^{q+2}}\|\nabla u\|_p^{p-\theta(q+2)}\right]^{\frac{2}{(1-\theta)(q+2)}}
\\
&\geq\left[\frac{a}{\beta^{q+2}}\tilde{\kappa}^{p-\theta(q+2)}\right]^{\frac{2}{(1-\theta)(q+2)}}.
\end{align*}

For the upper bound of $\Lambda_s$, using the H\"{o}lder inequality and \eqref{cao16}, we have
\begin{align*}
\Lambda_s&=\sup\limits_{u\in\mathscr{N}^s}\{\left\|u\|_2^2+k\|\nabla u\|_2^2\right\}
\\
&\leq\sup\limits_{u\in\mathscr{N}_s}(1+k)|\Omega|^{\frac{p-2}{p}}\|\nabla u\|_p^2
\\
&\leq (1+k)|\Omega|^{\frac{p-2}{p}}\tilde{\kappa}^2.
\end{align*}
$\hfill\Box$

\section{The ground state solutions}

\qquad Invoked by the ideas in \cite{Qu1} and \cite{Han3}, we introduced the following boundary value problem
\begin{eqnarray}\label{1.2}
\left\{
    \begin{array}{llll}
    \displaystyle -M(\|\nabla u\|_p^p)\Delta_p u=|u|^{q-1}u\log|u|, &&x\in\Omega,
    \\
    \displaystyle u(x)=0, &&x\in\partial\Omega.
    \end{array}
\right.
\end{eqnarray}
And the solution set
\begin{align*}
\Phi&=\{\text{the solutions of \eqref{1.2}}\}
\\
&=\{u\in W_0^{1,p}(\Omega): J'(u)=0 \text{in} W^{-1,p'}(\Omega)\}
\\
&=\{u\in W_0^{1,p}(\Omega): \langle J'(u), \varphi \rangle=0,
\forall\varphi\in W_0^{1,p}(\Omega)\}.
\end{align*}
In this section, we find the existence of ground state solutions to \eqref{1.2}
and establish some convergence relation between the global solutions of \eqref{1.1}
and the ground state solutions of \eqref{1.2}.

\begin{theorem}\label{th1.7}
There exists a function $u^*(x)\in\mathscr{N}$ such that
\\
{\rm(i)} $J(u^*)=\inf\limits_{u\in\mathscr{N}}J(u)=d$.
\\
{\rm(ii)} $u^*(x)$ is the ground state solution of \eqref{1.2},
namely $u^*\in\Phi\setminus\{0\}$ and $J(u^*)=\inf\limits_{u\in\Phi\setminus\{0\}}J(u)$.
\end{theorem}
{\bf Proof}  (i) From \eqref{3} and the definitions of $d$, $\mathscr{N}$, $J(u)$,
it is clear that there exists a minimizing sequence
$\{u_k\}_{k=1}^{\infty}\in\mathscr{N}$ such that
\begin{align*}
\lim\limits_{k\rightarrow \infty}J(u_k)=d
=\inf\limits_{u\in\mathscr{N}}\left\{\frac{a(q+1-p)}{p(q+1)}\|\nabla u\|_p^p
+\frac{b(q+1-2p)}{2p(q+1)}\|\nabla u\|_p^{2p}
+\frac{1}{(q+1)^2}\|u\|_{q+1}^{q+1}\right\},
\end{align*}
which induce that $\{u_k\}_{k=1}^{\infty}$ is bounded in $W_0^{1,p}(\Omega)$, namely
$$
\|u_k\|_{W_0^{1,p}}\leq M, \quad k=1,2,...
$$
with a positive constant $M$ independent of $k$. Due to the assumption
$1<2p-1<q<p^*-2$, $W_0^{1,p}(\Omega)\hookrightarrow L^{q+1}(\Omega)$
and $W_0^{1,p}(\Omega)\hookrightarrow L^{q+2}(\Omega)$ compactly.
Thus there exist a subsequence of $\{u_k\}_{k=1}^{\infty}$ (still denoted by $\{u_k\}_{k=1}^{\infty}$)
and $u^*\in W_0^{1,p}(\Omega)$, such that
\begin{center}
$u_k\rightharpoonup u^*$ weakly in $W_0^{1,p}(\Omega)$ as $k\to\infty$,
\\
$u_k\rightarrow u^*$ strongly in $L^{q+1}(\Omega)$ and $L^{q+2}(\Omega)$ as $k\to\infty$.
\end{center}
Next we prove that $\lim\limits_{{k\rightarrow+\infty}}\int_\Omega{|u_k|^{q+1}\log|u_k|}dx
=\int_\Omega{|u^*|^{q+1}\log|u^*|}dx$. We only need to prove that
\begin{align}\label{1.71}
\int_\Omega{\left||u_k|^{q+1}\log|u_k|-|u^*|^{q+1}\log|u^*|\right|}dx
\rightarrow0(k\rightarrow+\infty),
\end{align}
that is $\int_\Omega{\left((q+1)|\theta u_k
+(1-\theta)u^*|^q\log|\theta u_k+(1-\theta)u^*|
+|\theta u_k+(1-\theta)u^*|^q\right)\cdot\left|u_k-u^*\right|}dx
\rightarrow0(k\rightarrow+\infty)$, where $\theta\in(0,1)$.
Divide area $\Omega$ into two parts $\Omega=\Omega_1+\Omega_2$,
where $|\theta u_k+(1-\theta)u^*|\leq1$ in $\Omega_1$ and
$|\theta u_k+(1-\theta)u^*|>1$ in $\Omega_2$.

Case 1: $x\in\Omega_1$, then by the equality
$\inf{\left\{x^q\log x, x\in(0,1)\right\}}=-\frac{1}{eq}$, the H\"{o}lder inequality
and the boundedness of $u^*$, $u_k$, we have
\begin{align*}
&(q+1)\int_{\Omega_{1}}
{\left(|\theta u_k+(1-\theta)u^*|^q\log|\theta u_k+(1-\theta)u^*|\right)
\cdot\left|u_k-u^*\right|}dx
\\
&\leq\frac{q+1}{e}\int_{\Omega_{1}}
{\left(\frac{1}{eq}\left|u_k-u^*\right|\right)}dx
\\
&\leq\frac{(q+1)|\Omega_1|^{\frac{q+1}{q}}}{q}
\left(\int_{\Omega_{1}}{\left|u_k-u^*\right|^{q+1}}dx\right)^{\frac{1}{q+1}}
\\
&\rightarrow 0 (k\rightarrow +\infty).
\end{align*}

Case 2: $x\in\Omega_{2}$, then according to the H\"{o}lder inequality
and the boundedness of $u^*$, $u_k$, we have
\begin{align*}
&(q+1)\int_{\Omega_{2}}
{\left(|\theta u_k+(1-\theta)u^*|^q\log|\theta u_k
+(1-\theta)u^*|\right)\cdot\left|u_k-u^*\right|}dx
\\
&\leq(q+1)\int_{\Omega_{2}}
{\left(|\theta u_k+(1-\theta)u^*|^{q+1}\cdot\left|u_k-u^*\right|\right)}dx
\\
&\leq(q+1)\left(\int_{\Omega_{2}}{\left|u_k-u^*\right|^{q+2}}dx\right)^{\frac{1}{q+2}}
\cdot\left(\|u_k\|_{q+2}^{q+1}+\|u^*\|_{q+2}^{q+1}\right)
\\
&\rightarrow 0 (k\rightarrow +\infty).
\end{align*}
By the H\"{o}lder inequality, we also have
$\int_\Omega{\left(\left|\theta u_k+(1-\theta)u^*\right|^q\cdot\left|u_k-u^*\right|\right)}dx
\rightarrow 0 (k\rightarrow+\infty)$.
Thus we have $\lim\limits_{{k\rightarrow+\infty}}\int_\Omega{|u_k|^{q+1}\log|u_k|}dx
=\int_\Omega{|u^*|^{q+1}\log|u^*|}dx$.
\\
Since $\{u_k\}_{k=1}^{\infty}\in\mathscr{N}$, then we have $I(u_k)=0$, namely
$$
a\|\nabla u_k\|_p^p+b\|\nabla u_k\|_p^{2p}=\int_\Omega{|u_k|^{q+1}\log|u_k|}dx,
$$
which together with the weakly lower semi-continuity of $\|\cdot\|_{W_0^{1,p}}$ imply
\begin{align}\label{cao1}
a\|\nabla u^*\|_p^p+b\|\nabla u^*\|_p^{2p}
&\leq\liminf\limits_{k\rightarrow\infty}(a\|\nabla u_k\|_p^p+b\|\nabla u_k\|_p^{2p})\nonumber
\\
&=\lim\limits_{{k\rightarrow\infty}}\int_\Omega{|u_k|^{q+1}\log|u_k|}dx
\\
&=\int_\Omega{|u^*|^{q+1}\log|u^*|}dx.\nonumber
\end{align}
It is obviously that $\|u^*\|_{W_0^{1,p}}\neq0$. In what follows,
we assert that $u^*\in\mathscr{N}$, namely $I(u^*)=0$.
In fact, if it is not true, then \eqref{cao1} tells us that
$$
a\|\nabla u^*\|_p^p+b\|\nabla u^*\|_p^{2p}<\int_\Omega{|u^*|^{q+1}\log|u^*|}dx.
$$
By Lemma \ref{lem1}, there exists a unique $0<\lambda^*<1$
such that $\lambda^*u^*\in\mathscr{N}$ and $J(\lambda^*u^*)\geq d$. However,
\begin{align*}
J(\lambda^*u^*)&<\frac{a(q+1-p)}{p(q+1)}\|\nabla u^*\|_p^p
+\frac{b(q+1-2p)}{2p(q+1)}\|\nabla u^*\|_p^{2p}
+\frac{1}{(q+1)^2}\|u^*\|_{q+1}^{q+1}
\\
&\leq\liminf\limits_{{k\rightarrow+\infty}}
\frac{a(q+1-p)}{p(q+1)}\|\nabla u_k\|_p^p
+\frac{b(q+1-2p)}{2p(q+1)}\|\nabla u_k\|_p^{2p}
+\frac{1}{(q+1)^2}\|u_k\|_{q+1}^{q+1}
\\
&\leq d,
\end{align*}
which is a contradiction.
Thus \eqref{cao1} indicates that
$$
a\|\nabla u^*\|_p^p+b\|\nabla u^*\|_p^{2p}
-\int_\Omega{|u^*|^{q+1}\log|u^*|}dx
=\lim\limits_{k\rightarrow\infty}a\|\nabla u_k\|_p^p
+b\|\nabla u_k\|_p^{2p}-\int_\Omega{|u_k|^{q+1}\log|u_k|}dx,
$$
which with the uniform convexity of $W_0^{1,p}(\Omega)$ and
the weak convergence of $u_k$ in $W_0^{1,p}(\Omega)$,
we get $u_k\rightarrow u^*$ strongly in $W_0^{1,p}(\Omega)$ \cite{BrezisBOOK}. Then
$$
J(u^*)={\frac{a}{p}}\|\nabla u^*\|_p^p+{\frac{b}{2p}}\|\nabla u^*\|_p^{2p}
-{\frac{1}{q+1}}\int_\Omega{|u^*|^{q+1}\log|u|}dx+\frac{1}{(q+1)^2}\|u^*\|_{q+1}^{q+1}
=\lim\limits_{k\rightarrow \infty}J(u_k)=d,
$$
which implies $J(u^*)=\inf\limits_{u\in\mathscr{N}}J(u)=d$.

(ii) By conclusion (i) and $u^*\in\mathscr{N}$, we have
$\langle J'(u^*), u^*\rangle=I(u^*)=0$ and $J(u^*)=d$
with $\|u^*\|_{W_0^{1,p}}\neq 0$. According to the Lagrange multiplier method,
there exists a constant $\mu\in\mathbb{R}$ such that
\begin{align*}
J'(u^*)-\mu I'(u^*)=0,
\end{align*}
which implies $\mu\langle I'(u^*), u^*\rangle=\langle J'(u^*), u^*\rangle=0$.
For any $\varphi\in W_0^{1,p}(\Omega)$, we can deduce that
\begin{align*}
\langle I'(u^*), \varphi\rangle=&\frac{d}{d\tau}I\left(u^*+\tau\varphi\right)\Big|_{\tau=0}
\\
=&ap\left(|\nabla u^*|^{p-2}\cdot\nabla u^*, \nabla\varphi\right)
+2b p\|\nabla u^*\|_p^p\left(|\nabla u^*|^{p-2}\cdot\nabla u^*, \nabla\varphi\right)
\\
&-(q+1)\left(|u^*|^{q-1}\cdot u^*\log|u^*|, \psi\right)-\left(|u^*|^{q-1}\cdot u^*, \psi\right).
\end{align*}
Choosing $\varphi=u^*$ in the above equality leads to
\begin{align*}
\langle I'(u^*), u^*\rangle
=ap\|\nabla u^*\|_p^p+2bp\|\nabla u^*\|_p^{2p}-(q+1)\int_\Omega{|u^*|^{q+1}\log|u^*|}dx
-\|u^*\|_{q+1}^{q+1},
\end{align*}
which together with $I(u^*)=0$ points that
\begin{align*}
\langle I'(u^*), u^*\rangle=a(p-q-1)\|\nabla u^*\|_p^p
+b(2p-q-1)\|\nabla u^*\|_p^{2p}-\|u^*\|_{q+1}^{q+1}<0.
\end{align*}
Consequently $\mu=0$ and $J'(u^*)=0$, which mean $u^*\in\Phi\setminus\{0\}$.
From $\Phi\setminus\{0\}\subset\mathscr{N}$ and the conclusion (i), we get
$$
J(u^*)=\inf\limits_{u\in\Phi\setminus\{0\}}J(u)=d,
$$
which means that $u^*(x)$ is the ground state solution of \eqref{1.2}.
$\hfill\Box$

\begin{theorem}\label{th1.8}
Let $u(x,t)$ be the global solution of \eqref{1.1}.
Then there exists a function $u^*(x)\in\Phi$ and an increasing sequence
$\{t_k\}_{k=1}^\infty$ with $t_k\rightarrow +\infty$
when $k\rightarrow \infty$, such that
$$
\lim\limits_{k\rightarrow +\infty}\left\|\nabla\left(u(\cdot,t_k)-u^*\right)\right\|_p=0.
$$
\end{theorem}
{\bf Proof}  Using \eqref{cao2}, we have $J(u)\leq J(u_0)$, $t\in[0,+\infty)$.
Without loss of generality, we may further assume $J(u)\geq0$, $t\in[0,+\infty)$.
Otherwise, if there exists a $t_0>0$ such that $J(u(t_0))<0$.
Then by \eqref{3}, we have $I(u(t_0))<0$, which means that $u$ will blows up in finite time.
This contradicts the assumption that $u$ is global. Thus
\begin{equation}\label{cao4}
0\leq J(u)\leq J(u_0), \quad t\in[0,+\infty).
\end{equation}
Since $J(u)$ is non-increasing with respect to $t$,
there exists a constant $C\in[0,J(u_0)]$ such that
$\lim\limits_{t\rightarrow+\infty}J(u(t))=C$.
Integrating \eqref{cao2} from $0$ to $\infty$, we get
\begin{align*}
\int_0^{+\infty}(\|u_t\|_2^2+\|\nabla u_t\|_2^2)dt=J(u_0)-C\leq J(u_0),
\end{align*}
which indicates that there is an increasing sequence $\{t_k\}_{k=1}^{\infty}$
with $t_k\to\infty$ as $k\to\infty$ such that
$$
\lim\limits_{k\rightarrow \infty}\|u_t(t_k)\|_2^2+\|\nabla u_t(t_k)\|_2^2=0.
$$
On the other hand, for any $\varphi\in W_0^{1,p}(\Omega)$, we have
\begin{align*}
\langle J'(u), \varphi\rangle&=\frac{d}{d\tau}J(u+\tau\varphi)\Big|_{\tau=0}
\\
&=a\left(|\nabla u|^{p-2}\cdot\nabla u, \nabla\varphi\right)
+b\|\nabla u\|_p^p\left(|\nabla u|^{p-2}\cdot\nabla u, \nabla\varphi\right)
-\left(|u|^{q-1}\cdot u\log|u|,\varphi\right)
\\
&=\left(-u_t+\Delta u_t, \varphi\right),
\end{align*}
which with the H\"{o}lder inequality and the Sobolev imbedding inequality imply that
\begin{align}\label{cao6}
\|J'(u(t_k))\|_{W^{-1,p'}}&
=\sup\limits_{\|\varphi\|_{W_0^{1,p}}\leq 1}
\left|\langle J'(u(t_k)), \varphi\rangle\right|\nonumber
\\
&\leq\sup\limits_{\|\varphi\|_{W_0^{1,p}}\leq 1}
\left[\|u_t(t_k)\|_2\|\varphi\|_2
+\|\nabla u_t(t_k)\|_2\|\nabla\varphi\|_2\right]
\\
&\leq\sup\limits_{\|\varphi\|_{W_0^{1,p}}\leq 1}
C\left[\|u_t(t_k)\|_2\|\nabla\varphi\|_p
+\|\nabla u_t(t_k)\|_2\|\nabla\varphi\|_p\right]\nonumber
\\
&\rightarrow 0,\qquad k\rightarrow \infty,\nonumber
\end{align}
where $C$ is a constant independent of $k$.
From the above inequality and \eqref{2}, there exists a $\sigma>0$ such that
\begin{align*}
\frac{1}{q+1}|I(u(t_k))|&=\frac{1}{q+1}\left|\langle J'(u(t_k)), u(t_k)\rangle\right|
\\
&\leq\frac{1}{q+1}\|J'(u(t_k))\|_{W^{-1,p'}}\|u(t_k)\|_{W_0^{1,p}}
\\
&\leq\sigma\|\nabla u(t_k)\|_p,
\end{align*}
which with \eqref{cao4} guarantees that
\begin{align*}
J(u_0)+\sigma\|\nabla u(t_k)\|_p&\geq J(u(t_k))-\frac{1}{q+1}I(u(t_k))
\\
&=\frac{a(q+1-p)}{p(q+1)}\|\nabla u(t_k)\|_p^p
+\frac{b(q+1-2p)}{2p(q+1)}\|\nabla u(t_k)\|_p^{2p}
+\frac{1}{(q+1)^2}\|u(t_k)\|_{q+1}^{q+1}.
\end{align*}
Thus we have
$$
\frac{a(q+1-p)}{p(q+1)}\|\nabla u(t_k)\|_p^p-\sigma\|\nabla u(t_k)\|_p-J(u_0)
\leq0.
$$
Let $h(f)=\frac{a(q+1-p)}{p(q+1)}f^p-\sigma f-J(u_0)$ for $f\geq 0$.
Then $h(0)=-J(u_0)<0$, $h(+\infty)=+\infty$, $h(f)$ is strictly increasing on
$\left[\left(\frac{\sigma(q+1)}{a(q+1-p)}\right)^{\frac{1}{p-1}}, +\infty\right)$
and strictly decreasing on
$\left[0, \left(\frac{\sigma(q+1)}{a(q+1-p)}\right)^{\frac{1}{p-1}}\right)$.
Therefore there exists a unique $\hat{\kappa}>0$ such that $h(\hat{\kappa})=0$, which means
\begin{equation}\label{cao5}
\|\nabla u(t_k)\|_p\leq\hat{\kappa}, k=1,2,....
\end{equation}
Hence there exist $u^*\in W_0^{1,p}(\Omega)$ and
an increasing subsequences of $\{t_k\}_{k=1}^\infty$,
still represented by $\{t_k\}_{k=1}^\infty$, such that
\begin{center}
$u(t_k)\rightharpoonup u^*$ weakly in $W_0^{1,p}(\Omega)$ as $k\to\infty$,
\\
$u(t_k)\rightarrow u^*$ strongly in $L^{q+1}(\Omega)$ and $L^{q+2}(\Omega)$ as $k\to\infty$.
\end{center}
In what follows, we claim that $u(t_k)\rightarrow u^*$ strongly in $W_0^{1,p}(\Omega)$.
It is obvious that
\begin{align*}
\langle J'(u(t_k)), u(t_k)-u^*\rangle=&\frac{d}{d\tau}J\left(u(t_k)
+\tau\left(u(t_k)-u^*\right)\right)\Big|_{\tau=0}
\\
=&a\left(|\nabla u(t_k)|^{p-2}\cdot\nabla u(t_k), \nabla(u(t_k)-u^*)\right)
\\
&+b\|\nabla u(t_k)\|_p^p\cdot\left(|\nabla u(t_k)|^{p-2}\cdot\nabla u(t_k), \nabla(u(t_k)-u^*)\right)
\\
&-\left(|u(t_k)|^{q-1}\cdot u(t_k)\log|u(t_k)|,u(t_k)-u^*\right)
\end{align*}
and
\begin{align*}
\langle J'(u^*), u(t_k)-u^*\rangle
=&\frac{d}{d\tau}J\left(u^*+\tau\left(u(t_k)-u^*\right)\right)\Big|_{\tau=0}
\\
=&a\left(|\nabla u^*|^{p-2}\cdot\nabla u^*,\nabla(u(t_k)-u^*)\right)
\\
&+b\|\nabla u^*\|_p^p\cdot\left(|\nabla u^*|^{p-2}\cdot\nabla u^*, \nabla(u(t_k)-u^*)\right)
\\
&-\left(|u^*|^{q-1}\cdot u^*\log|u^*|, u(t_k)-u^*\right).
\end{align*}
It can be derived from \eqref{cao6} and \eqref{cao5} that
\begin{align}\label{cao7}
|\langle J'(u(t_k)),u(t_k)-u^*\rangle|
&\leq\|J'(u(t_k))\|_{W^{-1,p'}}\cdot\left(\|u(t_k)\|_{W_0^{1,p}}+\|u^*\|_{W_0^{1,p}}\right)\nonumber
\\
&\leq\left(\hat{\kappa}+\|u^*\|_{W_0^{1,p}}\right)\cdot\|J'(u(t_k))\|_{W^{-1,p'}}
\\
&\rightarrow 0,\qquad\hbox{as}\quad k\rightarrow \infty.\nonumber
\end{align}
Since $u(t_k)\rightharpoonup u^*$ weakly in $W_0^{1,p}(\Omega)$, we also have
\begin{equation}\label{cao8}
\langle J'(u^*), u(t_k)-u^*\rangle\rightarrow 0, \qquad\hbox{as}\quad k\rightarrow \infty.
\end{equation}
Combining \eqref{cao7} and \eqref{cao8} leads to
\begin{align*}
\langle J'(u(t_k))-J'(u^*), u(t_k)-u^*\rangle
=&a\left(|\nabla u(t_k)|^{p-2}\cdot\nabla u(t_k)-|\nabla u^*|^{p-2}\cdot\nabla u^*, \nabla(u(t_k)-u^*)\right)
\\
&+b\|\nabla u(t_k)\|_p^p\left(|\nabla u(t_k)|^{p-2}\cdot\nabla u(t_k), \nabla(u(t_k)-u^*)\right)
\\
&-b\|\nabla u^*\|_p^p\left(|\nabla u^*|^{p-2}\cdot\nabla u^*, \nabla(u(t_k)-u^*)\right)
\\
&-\left(|u(t_k)|^{q-1}\cdot u(t_k)\log|u(t_k)|-|u^*|^{q-1}\cdot u^*\log|u^*|, u(t_k)-u^*\right)
\\
=&M(\|\nabla u(t_k)\|_p^p)\left(|\nabla u(t_k)|^{p-2}\cdot\nabla u(t_k), \nabla(u(t_k)-u^*)\right)
\\
&-M(\|\nabla u^*\|_p^p)\left(|\nabla u^*|^{p-2}\cdot\nabla u^*, \nabla(u(t_k)-u^*)\right)
\\
&-\left(|u(t_k)|^{q-1}\cdot u(t_k)\log|u(t_k)|-|u^*|^{q-1}\cdot u^*\log|u^*|, u(t_k)-u^*\right)
\\
\rightarrow& 0, \qquad\hbox{as}\quad k\rightarrow +\infty.
\end{align*}
Similar to the proof method of \eqref{1.71}, we have
\begin{align*}
\left|\left(|u(t_k)|^{q-1}\cdot u(t_k)\log|u(t_k)|-|u^*|^{q-1}\cdot u^*\log|u^*|, u(t_k)-u^*\right)\right|
\rightarrow 0, \qquad\hbox{as}\quad k\rightarrow \infty.
\end{align*}
Using $u(t_k)\rightharpoonup u^*$ weakly in $W_0^{1,p}(\Omega)$ again, we can get
$$
M(\|\nabla u^*\|_p^p)\left(|\nabla u^*|^{p-2}\cdot\nabla u^*, \nabla(u(t_k)-u^*)\right)
\rightarrow 0, \qquad\hbox{as}\quad k\rightarrow \infty.
$$
Then we can deduce that
\begin{align}\label{cao9}
M(\|\nabla u(t_k)\|_p^p)\left(|\nabla u(t_k)|^{p-2}\cdot\nabla u(t_k), \nabla(u(t_k)-u^*)\right)
\rightarrow 0, \qquad\hbox{as}\quad k\rightarrow \infty.
\end{align}
From the weak convergence of $\{u(t_k)\}_{k=1}^{\infty}$, we have
\begin{align*}
\|\nabla u^*\|_p&\leq\liminf\limits_{k\rightarrow \infty}\|\nabla u(t_k)\|_p.
\end{align*}
If there exists a subsequence of $\{t_k\}_{k=1}^{\infty}$, denoted still by itself,
such that $\lim\limits_{k\rightarrow \infty}\|\nabla u(t_k)\|_p=0$,
then $\|\nabla u^*\|_p=0$, namely $u^*=0$, the conclusion is proved.
If not and $\lim\limits_{k\rightarrow \infty}\|\nabla (u(t_k)-u^*)\|_p\neq 0$,
then there exist $\epsilon>0$ and a subsequence of $\{t_k\}_{k=1}^{\infty}$,
still represented by itself, such that
\begin{align*}
\|\nabla u^*\|_p<\|\nabla u(t_k)\|_p-\epsilon.
\end{align*}
Then we have
\begin{align*}
&M(\|\nabla u(t_k)\|_p^p)\left(|\nabla u(t_k)|^{p-2}\cdot\nabla u(t_k), \nabla(u(t_k)-u^*)\right)
\\
=&M(\|\nabla u^*\|_p^p)\int_{\Omega}\left(|\nabla u(t_k)|^p-|\nabla u(t_k)|^{p-2}\cdot\nabla u(t_k)\cdot\nabla u^*\right)dx
\\
\geq&M(\|\nabla u^*\|_p^p)\left(\|\nabla u(t_k)\|_p^p-\int_{\Omega}{|\nabla u(t_k)|^{p-1}\cdot|\nabla u^*|}dx\right)
\\
\geq&M(\|\nabla u^*\|_p^p)\left(\|\nabla u(t_k)\|_p^p-\|\nabla u(t_k)\|_p^{p-1}\cdot\|\nabla u^*\|_p\right)
\\
=&M(\|\nabla u^*\|_p^p)\|\nabla u(t_k)\|_p^{p-1}\left(\|\nabla u(t_k)\|_p-\|\nabla u^*\|_p\right)
\\
\geq&M(\|\nabla u^*\|_p^p)\cdot\|\nabla u(t_k)\|_p^{p-1}\cdot\epsilon,
\end{align*}
which contradicts \eqref{cao9}.
So we conclude that $u(t_k)\rightarrow u^*$ strongly in $W_0^{1,p}(\Omega)$.
Consequently, $J'(u^*)=\lim\limits_{k\rightarrow \infty}J'(u(t_k))=0$, namely $u^*\in\Phi$.
$\hfill\Box$

\end{document}